\documentclass[11pt]{article}
\usepackage{amssymb,amsthm,amsmath} 
\usepackage[left=2.5cm,right=2.5cm,top=2.5cm,bottom=2.5cm]{geometry}
\usepackage{amsmath}
\usepackage{soul} \usepackage{url}

\usepackage{lineno,xcolor}

\usepackage{amssymb}
\def\mymedskip{\vskip\medskipamount}
\def\mymedbreak{\par \ifdim\lastskip<\medskipamount
  \removelastskip \penalty-100 \mymedskip \fi}
\def\myaftermedspace{\par \ifdim\lastskip<\medskipamount
  \removelastskip \penalty55\mymedskip\fi}
\newcommand{\eop}{{\unskip\nobreak\hfil\penalty50
          \hskip2em\hbox{}\nobreak\hfil$\Box$
          \parfillskip=0pt \finalhyphendemerits=0 \par}}
{\mymedbreak{\noindent\bf Proof of Theorem #1:\enspace}}{\eop\myaftermedspace}
\newenvironment{remark}%
{\mymedbreak\noindent{\bf Remark:}%
\enspace\rm}%
{\myaftermedspace}
\theoremstyle{definition}
\newtheorem{teor}{Theorem}[section]
\newtheorem{defi}[teor]{Definition}
\newtheorem{fact}[teor]{Fact}
\newtheorem{problem}{Problem}
\newtheorem{exercise}{Exercise}
\newtheorem{examp}[teor]{Example}
\newtheorem{lem}[teor]{Lemma}
\newtheorem{cor}[teor]{Corollary}

\newtheorem{prop}[teor]{Proposition}
\newtheorem{rem}[teor]{Remark}
\newcommand{\beq}{\begin{equation}}
\newcommand{\eeq}{\end{equation}}
\newcommand{\beql}[1]{\begin{equation} \label{#1}}
\newcommand{\eeql}{\end{equation}}
\newcommand{\beqa}{\begin{eqnarray*}}
\newcommand{\eeqa}{\end{eqnarray*}}
\newcommand{\beqal}[1]{\begin{eqnarray} \label{#1}}
\newcommand{\eeqal}{\end{eqnarray}}
\newcommand{\beqan}{\begin{eqnarray}}
\newcommand{\eeqan}{\end{eqnarray}}
\newcommand{\bpf}{\begin{proof}}
\newcommand{\epf}{\end{proof}}
\newcommand{\ben}{\begin{enumerate}}
\newcommand{\een}{\end{enumerate}}
\newcommand{\bit}{\begin{itemize}}
\newcommand{\eit}{\end{itemize}}

\newcommand{\cB}{{\cal B}}

\newcommand{\cE}{{\cal E}}

\newcommand{\cR}{{\cal R}}

\newcommand{\cZ}{{\cal Z}}

\newcommand{\bZ}{{\bf Z}}

\newcommand{\GF}{{\rm GF}}

\newcommand{\gs}{\sigma}

\newcommand{\ga}{\alpha}
\newcommand{\gb}{\beta}
\newcommand{\gc}{\gamma}
\newcommand{\gd}{\delta}
\newcommand{\gD}{\Delta}
\newcommand{\gre}{\epsilon}
\newcommand{\gl}{\lambda}

\newcommand{\bab}{\begin{abstract}}
\newcommand{\eab}{\end{abstract}}
\newcommand{\bke}{\begin{keywords}}
\newcommand{\eke}{\end{keywords}}

\newcommand{\btm}[1]{\begin{teor} \label{#1}}
\newcommand{\etm}{\end{teor}}
\newcommand{\btmn}[2]{\begin{teor}[#1] \label{#2}}
\newcommand{\etmn}{\end{teor}}
\newcommand{\ble}[1]{\begin{lem} \label{#1}}
\newcommand{\ele}{\end{lem}}
\newcommand{\bLe}[1]{\begin{Lemma} \label{#1}}
\newcommand{\eLe}{\end{Lemma}}
\newcommand{\bpn}[1]{\begin{prop} \label{#1}}
\newcommand{\epn}{\end{prop}}
\newcommand{\bex}[1]{\begin{examp} \label{#1}}
\newcommand{\eex}{\end{examp}}
\newcommand{\bde}[1]{\begin{defi} \label{#1}}
\newcommand{\ede}{\end{defi}}
\newcommand{\bco}[1]{\begin{cor} \label{#1}}
\newcommand{\eco}{\end{cor}}
\newcommand{\bcorn}[2]{\begin{cor}[#1] \label{#1}}
\newcommand{\ecorn}{\end{cor}}
\newcommand{\bcon}[1]{\begin{conjecture} \label{#1}}
\newcommand{\econ}{\end{conjecture}}
\newcommand{\bfa}[1]{\begin{fact} \label{#1}}
\newcommand{\efa}{\end{fact}}
\newcommand{\bpr}[1]{\begin{problem} \label{#1}}
\newcommand{\epr}{\end{problem}}
\newcommand{\bexer}[1]{\begin{exercise} \label{#1}}
\newcommand{\eexer}{\end{exercise}}
\newcommand{\bre}[1]{\begin{rem} \label{#1}}
\newcommand{\ere}{\end{rem}}

\newcommand{\bbF}{\mathbb{F}}

\newcommand{\bbP}{\mathbb{P}}
\newcommand{\bbQ}{\mathbb{Q}}

\newcommand{\bbZ}{\mathbb{Z}}

\newcommand{\diag}{{\rm diag}}

\newcounter{question_number}
\newenvironment{question}{\addtocounter{question_number}{1}\noindent{\bf Question \arabic{question_number}}}{\myaftermedspace}

\newenvironment{solution}{\noindent {\bf Solution:} \enspace}{\eop\myaftermedspace}

\newenvironment{hint}{\noindent {\bf Hint:} \enspace}{\eop\myaftermedspace}

\newenvironment{multisolution}[1]{\noindent {\bf Solution #1:} \enspace}{\eop\myaftermedspace}

\newcommand{\bqu}{\begin{question}}
\newcommand{\equ}{\end{question}}
\newcommand{\bs}{\begin{solution}}
\newcommand{\es}{\end{solution}}
\newcommand{\bh}{\begin{hint}}
\newcommand{\eh}{\end{hint}}
\newcommand{\bms}[1]{\begin{multisolution}{#1}}
\newcommand{\ems}{\end{multisolution}}

        \DeclareMathOperator{\GL}{GL}

 \newcommand{\FF}{\mathbb{F}}

 \newcommand{\te}{\tilde{e}}

 \newcommand{\ord}{{\rm ord}}
\newcommand{\gS}{\Sigma} \newcommand{\im}{{\rm Im}}
 \renewcommand{\ker}{{\rm Ker}}
\newcommand{\cV}{{\cal V}} 
\newcommand{\cVs}{\cV^*} 

\begin{document} 
\begin{titlepage} 
\title{Sandpile groups of
generalized de Bruijn  and Kautz graphs and 
circulant matrices over finite fields} 
\date{\today} 
\author{%
Swee Hong Chan\thanks{
Division of Mathematical Sciences, School of Physical and
Mathematical Sciences, Nanyang Technological University, Singapore, email:
\{sweehong, henk.hollmann\}@ntu.edu.sg}
\and Henk D.L.~ Hollmann\footnotemark[1]
\and Dmitrii V.~ Pasechnik\thanks{
Department of Computer Science, University of Oxford, UK, email:
dmitrii.pasechnik@cs.ox.ac.uk (corresponding author)}
}
\maketitle \begin{abstract}
A maximal minor $M$ of the Laplacian of an $n$-vertex Eulerian digraph
$\Gamma$ gives rise to
a finite group $\bbZ^{n-1}/\bbZ^{n-1}M$ known as the 
 {\em sandpile} (or {\em critical}) {\em group} $S(\Gamma)$ of $\Gamma$. 
We determine $S(\Gamma)$ of the {\em generalized
de Bruijn graphs} $\Gamma=\mathrm{DB}(n,d)$ with vertices $0,\dots,n-1$
and arcs $(i,di+k)$ for $0\leq i\leq n-1$ and $0\leq k\leq d-1$,
and closely related {\em generalized Kautz graphs},
extending and completing earlier results
for the classical de Bruijn and Kautz graphs. 

Moreover, for a prime $p$ and an $n$-cycle permutation matrix 
$X\in\mathrm{GL}_n(p)$ we show that $S(\mathrm{DB}(n,p))$ is isomorphic to
the quotient by $\langle X\rangle$
of the centraliser of $X$ in $\mathrm{PGL}_n(p)$. 
This offers an explanation for the coincidence of
numerical data in sequences A027362 and A003473 of the OEIS, and 
allows one to speculate upon a possibility to construct normal
bases in the finite field $\mathbb{F}_{p^n}$ from spanning trees
in~$\mathrm{DB}(n,p)$.
\end{abstract} 

\begin{center}
{\it Dedicated to the memory of \'{A}kos Seress.}
\end{center}

\end{titlepage}
%


\section{\label{sect:int}Introduction}
The {\em critical group\/} $S(G)$ of a directed graph~$G$ is an
abelian group obtained from the Laplacian matrix $\gD$ of~$G$. It
carries the same information as the Smith Normal Form (SNF) of~$\gD$.
(For precise definitions of these and other terms, we refer to the
next section.) The {\em sandpile group\/} $S(G,v)$ of a
directed graph or {\em digraph\/}~$G$  at a vertex $v$ is an abelian
group obtained from the reduced Laplacian $\gD_v$ of~$G$; by the
Matrix Tree Theorem \cite{AaDB-tree}, its order is equal to 
the number of directed 
trees rooted at~$v$, see for example \cite{LeLI-sand}.
If the graph~$G$ is Eulerian, then the sandpile group does not depend
on the vertex~$v$ and is equal to the critical group $S(G)$ of~$G$
\cite{HLMPPW-chip}.  Most of the literature on sandpile groups is
concerned with
{\em undirected\/} graphs, 
which can be considered as a special case, namely for directed graphs
that are obtained by replacing each undirected edge  in a graph by a
pair of directed edges oriented in opposite directions. 

The critical group has been studied in 
other contexts under variuos other names, such as  group of components
(in arithmetic geometry), 
Jacobian group and Picard group (for algebraic curves), and Smith
group (for matrices). For more details and background, see, for
example, \cite{Lo-fin, Ru-phd, AlVa-cone} for the undirected case and
\cite{ HLMPPW-chip,  wag-crit} for the directed case. 

Critical groups have been determined for many families of (mostly
undirected) graphs. For some examples,  see \cite{ChHou-pc,
HWC-sqc,Ra-W3, RaWa-W4,CHW-mob,DaFiFr-dihe,SH-bra, Me-eq,Bi-chip}, and
the references in \cite{AlVa-cone}.  
Here, we determine the
critical group of the generalized de Bruijn graphs $\mathrm{DB}(n,d)$
and generalized Kautz graphs $\mathrm{Ktz}(n,d)$ (which are in fact
both directed graphs), thereby extending and completing the results
from \cite{LeLI-sand} for the binary de Bruijn graphs
$\mathrm{DB}(2^\ell,2)$ and Kautz graphs
$\mathrm{Ktz}((p-1)p^{\ell-1},p)$ for primes $p$, and
\cite{BiKi-knuth} for the classical de Bruijn graphs
$\mathrm{DB}(d^\ell,d)$ and Kautz graphs
$\mathrm{Ktz}((d-1)d^{\ell-1},d)$. Unlike the classical case, the
generalized versions are not necessarily  iterated line graphs, so to
obtain their critical groups, techniques different from these in
\cite{LeLI-sand} and \cite{BiKi-knuth} have to be applied. 

As set out  in \cite{DuPa-neck}, our original motivation for studying
the sandpile  groups of generalized de Bruijn graphs
was to explain their apparent relation to some other algebraic
objects, such as the groups $C(n,p)$ of invertible $n\times
n$-circulant matrices over $\mathbb{F}_{p}$ (mysterious numerical
coincidences of the OEIS entries A027362 and A003473  
were noted by the third author \cite{oeis:A027362}, 
while running extensive computer experiments using
Sage~\cite{sage,sagesandpiles}), and {\em normal bases} (cf. e.g.
\cite{MR1429394}) of finite fields $\mathbb{F}_{p^n}$, in the case
where $(n,p)=1$. The latter were noted to be closely related to
circulant matrices and to {\em necklaces\/} by Reutenauer
\cite[Sect.~7.6.2]{Reu-lie}, see also \cite{DuPa-neck} and the related
numeric data collected in~\cite{Arn-comp}.  Here we show that the
critical group of $\mathrm{DB}(n,p)$ is isomorphic to
$C(n,p)/\langle Q_n\rangle \times \bbF_p^*$, where $Q_n$ denotes the
permutation matrix of the $n$-cycle. Although we were not able
to construct an explicit bijection between the former and the latter,
we could speculate that potentially one might be able to design a new
deterministic way to construct normal bases of  $\mathbb{F}_{p^n}$
from spanning trees in~$\mathrm{DB}(n,p)$.
For more details and background on this connection with aperiodic
necklaces, we refer the interested reader to~\cite{DuPa-neck}.  Most
of the results for the case $d$ prime were first derived in
the undergraduate thesis  \cite{chan-db} by the first author, 
supervised by the third author. Results in this text were announced in
an extended abstract \cite{sandpiles_eurocomb13}.

\section{Preliminaries}\label{sect:prelims}
In this section, we introduce the necessary terminology and
background. First, in Section~\ref{ssect:scs},  we  discuss the Smith
Normal
Form and the Smith group, and the critical group and sandpile group of
a directed graph,  as well as the  relations between these notions.
Then in
Section~\ref{ssect:gdbk}, we define the generalized De Bruijn and
Kautz graphs. The group of invertible circulants  is defined in
Section~\ref{ssect:Icm}. Finally, in Section~\ref{ssect:sd}, we derive
expressions for the sandpile group of generalized De Bruijn and Kautz
graphs as embeddings in a group that we refer to as sand dune
group.

\subsection{\label{ssect:scs}Smith group, Critical Group, and Sandpile
Group}
Let $M$ be an rank $r$ integer $m\times n$ matrix. There exist 
positive integers $s_1, \ldots, s_r$ 
with $s_i|s_{i+1}$ for $i=1, \ldots, r$ 
and unimodular 
matrices $P$ and $Q$
such that $PMQ=D=\diag(s_1, \ldots, s_r, 0, \ldots, 0)$.
The diagonal matrix $D$ is called the {\em Smith Normal Form\/} (SNF)
of $M$, and the numbers $s_1, \ldots, s_r$ are the 
nonzero {\em invariant factors\/} of~$M$.  
The SNF, and hence the invariant factors, are  {\em uniquely\/}
determined by the matrix~$M$.  For background on the SNF and invariant
factors, see, e.g.,~\cite{Ne-im}.
The {\em Smith group\/}~\cite{Ru-phd} of $M$ is 
$\Gamma(M)=\bbZ^n/\bbZ^m M$;
the submodule
$\overline{\Gamma}(M)=\bbZ^n/\bbQ^mM\cap \bbZ^n$ of~$\Gamma(M)$ is the 
torsion subgroup of $\Gamma(M)$. 
Indeed, if $M$ has rank~$r$, then
$\Gamma(M) = \bbZ^{n-r} \oplus \overline{\Gamma}(M)$
with 
$\overline{\Gamma}(M) = \oplus_{i=1}^r \bbZ_{s_i}$, where $s_1,
\ldots, s_r$ are the nonzero invariant factors of~$M$.
See \cite{Ru-phd} for further details and proofs.

Let  $G=(V,E)$ be a finite directed graph with
with vertex set $V$ and edge set $E$ (we allow loops and multiple
edges). The {\em adjacency matrix\/} of~$G$ is the $|V|\times |V|$
matrix $A=(A_{v,w})$, with rows and columns indexed by $V$, where the
entry $A_{v,w}$ in position $(v,w)$  is the number of edges from $v$
to~$w$.  The indegree $d^-(v)$ and outdegree $d^+(v)$ is the number of
edges ending or starting in vertex~$v$, respectively.
The {\em Laplacian\/} of~$G$ is the matrix $\gD=D-A$, where $D$ is
diagonal with $D_{v,v}=d^+_v$.  The {\em critical group\/} $S(G)$
of~$G$ is the torsion subgroup $\overline{\Gamma}(\gD)$ of the Smith group
of the Laplacian~$\gD$ of~$G$.  The {\em sandpile group\/} $S(G,v)$
of~$G$ at a vertex~$v$ is the torsion subgroup of the Smith group of the 
{\em reduced Laplacian\/} $\gD_v$, obtained from $\gD$ by deleting the
row and the column of $\gD$ indexed by~$v$. Note that by the Matrix
Tree Theorem for directed graphs, the order of $S(G,v)$ equals the
number of directed spanning trees rooted at~$v$.  The directed graph
$G$ is called  {\em Eulerian\/} if $d^+(v)=d^-(v)$ for every
vertex~$v$. 
In that case, $S(G,v)$ does not depend on the vertex~$v$ and is equal
to the critical group $S(G)$ of~$G$, essentially because in that case,
not only the columns, but also the rows of the Laplacian $\gD$ of~$G$
sum to zero; for a detailed proof, see \cite{HLMPPW-chip} (note that
 the proof as given there does not use the assumption that the
directed graph is connected). For more details on sandpile groups and
the critical group of a directed graph, we refer for example to
\cite{HLMPPW-chip} or \cite{wag-crit}. 

\subsection{\label{ssect:gdbk}Generalized de Bruijn and Kautz graphs}
Generalized de Bruijn graphs \cite{DuHw-GdB} and generalized Kautz
graphs  \cite{DCH-dBKgen} are known to have a relatively small
diameter and attractive connectivity properties, and have been studied
intensively due to their applications in interconnection networks. The
generalized Kautz graphs were first investigated  in
\cite{ImIt-diam-ieee}, \cite{ImIt-diam-net}, and are
also known as {\em Imase-Itoh digraphs\/}. Both classes of directed
graphs are Eulerian.

Let $n$ and $d$ be integers with $n\geq1$ and $d\geq0$.
The {\em generalized de Bruijn graph\/} $\mathrm{DB}(n,d)$ has vertex
set $\bbZ_n$, the set of integers modulo~$n$, and $d$ (directed) edges
$v\rightarrow dv+i$,  for $0\leq i\leq d-1$, for $0\leq v\leq n-1$. 
The {\em generalized Kautz graph\/} 
${\rm Ktz}(n,d)$ has vertex set $\bbZ_n$ and $d$ directed edges
$v\rightarrow -d(v+1) +i$, for $0\leq i\leq d-1$,  for 
$0\leq v\leq n-1$. (For both generalized de Bruijn and Kautz graphs, we
allow multiple edges when $d>n$.)  These directed graphs are important
special cases of the so-called {\em consecutive-$d$ digraphs\/}. Here
a  consecutive-$d$ digraph $G(d, n, q, r)$, defined for $q\in
\bbZ_n\setminus \{0\}$ and $r\in\bbZ_n$, has vertex set $\bbZ_n$ and
directed edges $v\rightarrow qv+r +i$ for $0\leq i\leq d-1$ and
$0\leq v\leq n-1$. Note that the generalized de Bruijn and Kautz graphs
are the cases $q=d, r=0$ and $q=-d, r=-d$, respectively. 
It is easily verified that both $\mathrm{DB}(n,d)$ and
$\mathrm{Ktz}(n,d)$ are indeed Eulerian for all integers $n\geq1$ and
$d\geq0$. 
%
 It is easily seen that for $d=0$ or $1$, the critical groups for both
the de Bruijn graph and the Kautz graph are trivial.
%


 
\subsection{\label{ssect:Icm}The group of invertible circulant
matrices}
Let $q=p^r$ be a prime power. An $n\times n$ circulant matrix over a
finite field~$\bbF_q$ is a matrix~$C$ of the form 
\beql{Lcirc1def} 
C=
\left( \begin{array}{cccc} c_0& c_{n-1} & \cdots & c_1 \\ c_1& c_0 &
\ddots & \vdots  \\ \vdots & \ddots& \ddots  &c_{n-1}\\ c_{n-1}&
\cdots &c_1 & c_0 \end{array} 
\right), 
\eeql 
where $c_0, \ldots,
c_{n-1}\in \bbF_q$.  With $C$ as in (\ref{Lcirc1def})
we associate the polynomial $c_C(x):=c(x)=c_0+c_1x+\cdots +c_{n-1}x^{n-1}\in\bbF_q[x]$.
Let $X$ denote the matrix of the multiplication by $x$ map on 
the ring $\cR:=\bbF_q[x]/(x^n-1)$ with respect to
the basis $1, x, \ldots, x^{n-1}$.  Every $C$ 
in~(\ref{Lcirc1def})
can be written as $C=\sum_{i=0}^{n-1} c_iX^i$, where $c(x)=c_C(x)$.
Note that $c_X(x)=x$, and the map $x\mapsto X$ induces an
isomorphism between $\cR$ and the algebra $\bbF_q[X]$ of circular matrices. 

The units of $\bbF_q[X]$ form a commutative group
under multiplication; we denote this group by~$C(n,q)$. 
Under the isomorphism
induced by $x\mapsto X$, the group  $C(n,q)$ corresponds to
\[\cR^*=\{c(x)\in \bbF_q[x], \deg c<n \mid (c(x), x^n-1)=1\}.\]
Indeed, by the extended Euclidean algorithm, $(c(x), x^n-1)=1$ if and
only if there exist $u,v\in\bbF_q[x]$ such that $cu=1-(x^n-1)v$, i.e.
$u=c^{-1}\in \cR$. On the other hand, $c(X)u(X)=I-(X^n-I)v(X)=I$, 
i.e. $u(X)=c(X)^{-1}\in C(n,q)$.

Note that $C(n,q)$ contains a subgroup isomorphic to
$\mathbb{Z}_{q-1}\oplus \mathbb{Z}_n$, namely the direct product of
the group of scalar matrices $F_q^* I:=\{\lambda I\mid
\lambda\in\mathbb{F}_q^*\}$ and the cyclic subgroup $\langle
X\rangle$ generated by $X$. 
Each $C\in\bbF_q[X]$ has the all-ones vector
$\mathbf{1}:=(1,\dots,1)^\top$ as an eigenvector. Thus
$C'(n,q):=\{C\in C(n,q)\mid C\mathbf{1}=\mathbf{1}\}\leq C(n,q)$, 
and we have the following direct product decomposition.
\begin{equation}\label{Cnpdec} 
C(n,q)=C'(n,q)\times F_q^* I,
\end{equation} 
where, as usual, $\bbF_q^*=\bbF_q\setminus\{0\}$. In view of \eqref{Cnpdec}
one has $C'(n,q)\leq\mathrm{PLG}_n(q)$. Note that 
\beql{LECp}  
C'(n,q) \cong \{a(x)\in  \bbF_q[x], \deg a<n \mid
\mbox{$(a(x),x^n-1)=1$ and $a(1)=1$}\}.
\eeql  
We also note that although $\langle X\rangle\leq C'(n,q)$, it is  not
necessarily a {\em direct summand\/} of $C'(n,q)$ (for example, take
$n=2$ and $q=4$: then $C'(n,q)\cong\bbZ_4$ and $\langle X\rangle\cong\bbZ_2$).

\begin{remark}
Note that $C(n,q)$ is the centralizer of $X$ in $\GL_n(q)$.
It also can be viewed as the group of units of the group algebra
$\bbF_q[\langle X\rangle]$.
\end{remark}

\subsection{\label{ssect:sd}The sandpile and sand dune groups}
We determine the sandpile group $S_{\rm DB}(n,d)$ 
(resp. $S_{\rm Ktz}(n,d)$)
of a generalized de Bruijn (resp.  Kautz graph) on $n$ vertices 
by embedding this group as a subgroup of
index~$n$ in a group that, for lack of a better name, we refer to
as the {\em sand dune\/} group of the corresponding directed graph.
This embedding method can in fact be applied to the much wider class
of {\em consecutive-$d$ digraphs\/}~\cite{DHP-cons}, and the idea may
have other applications as well.

Let us 
represent the elements of $\bbZ_n$ by $1, x, \ldots x^{n-1}$,
and think of elements in the group algebra $\bbQ[{\bbZ_n}]$ as
Laurent polynomials modulo $x^n-1$, that is, identify
$\bbQ[{\bbZ_n}]$ with $\bbQ[x,x^{-1}] \bmod x^n-1$, and its subring
$\bbZ[{\bbZ_n}]$ with $\bbZ[x,x^{-1}] \bmod x^n-1$.  Furthermore, we
identify a vector $c=(c_0,\ldots, c_{n-1})$ in $\bbQ^n$ with its
associated polynomial $c(x)=c_0+\cdots +c_{n-1}x^{n-1}$ in
$\bbQ[x,x^{-1}] \bmod x^n-1$; note that this association is in fact an
isomorphism between $\bbQ^n$ and $\bbQ[x,x^{-1}] \bmod x^n-1$,
considered as vector spaces over $\bbQ$. The advantage of this
identification is that we now also have a multiplication available.
Given a collection of vectors $V=\{v_i \mid i\in I\}$ in $\bbZ^n$, we
denote the $\bbZ$-span of the associated polynomials by $\langle
v_i(x) \mid i\in I\rangle_{\bbZ}$.

We now derive a useful description of both $S_{\rm DB}(n,d)$ and
$S_{\rm Ktz}(n,d)$ using  the Smith group $\Gamma(\gD)$ of the
Laplacian $\gD$ of  these digraphs as defined in
Section~\ref{ssect:scs}.  To this end, for every integer $d$ we define
the Laurent polynomial $f^{(n,d)}(x)$ in $\bbZ[x,x^{-1}]\bmod x^n-1$
as 
\beql{LEfv}  f_v^{(n,d)}(x)= dx^v -x^{dv} (x^d-1)/x-1).
\eeql 
For $d\geq0$, we have
\[ f_v^{(n,d)}(x)= dx^v -x^{dv} \sum_{i=0}^{d-1} x^i =
\sum_{w\in \bbZ_n} \gD_{v,w}x^w  
\] 
with $\gD=\gD_{\rm DB}$ the Laplacian  of $\mathrm{DB}(n,d)$, while for
$d<0$, we have
\[ f_v^{(n,d)}(x)= -|d|x^v -x^{-|d|v}
  (x^{-|d|}-1)/(x-1)= -|d|x^v  + x^{-|d|(v+1)}
 \sum_{i=0}^{|d|-1} x^i =-\sum_{w\in \bbZ_n} \gD_{v,w}x^w  
\] 
with
$\gD=\gD_{\rm Ktz}$ now the Laplacian  of $\mathrm{Ktz}(n,|d|)$.
In what follows, we simply write $\gD$ to denote the
Laplacian of either $\mathrm{DB}(n,d)$  or $\mathrm{Ktz}(n,d)$.
For later use, we also define 
\[ g_v^{(n,d)}(x)= (x-1)f^{(n,d)}_v(x) =
dx^v (x-1) -x^{dv}(x^d-1)\] 
for every $0\leq v<n$.  For the
remainder of this paper, we let 
\begin{equation}\label{def:e_v}
e_v^{(n)}=x^v-1, \qquad \gre_v^{(n,d)}=de^{(n)}_v-e^{(n)}_{dv}, \quad
\text{for every $0\leq v\leq n-1$};
\end{equation}
note that $e_0^{(n)}=\gre_0^{(n,d)}=0$. We simply write $f_v(x)$,
$g_v(x)$, $e_v$ and $\gre_v$ if the intended values for $n$ and $d$
are evident from the context. Finally, we set 
\[\cZ_n =\langle e_v\mid
1\leq v\leq n-1\rangle_\bbZ, \qquad \cE_{n,d}=\langle
\gre_v\mid 1\leq v\leq n-1\rangle_\bbZ.\]
In the next lemma, we collect some simple facts.
\ble{LLbasic}\mbox{}\\
(1) We have that $\sum_{v=0}^{n-1}f_v(x)=0$ and $f_v(1)=0$;\\ 
(2) $\cZ_n$ consists of all polynomials $c(x)\in\bbZ[x]$
for which $\deg c\leq n-1$ and $c(1)=0$; \\
(3) We have that $\gre_v=g_0(x)+\cdots+g_{v-1}(x)$ and
$\cE_{n,d}=\langle  g_v(x)\mid 0<v<n \rangle_\bbZ$.
\ele \bpf 
(1) Since the columns of the Laplacian $\gD$ add
up to 0, we have that $f_v(1)=0$; moreover, since both
$\mathrm{DB}(n,d)$ and  $\mathrm{Ktz}(n,d)$ are Eulerian, in both
cases the  rows of $\gD$ also add up to 0, so that $\sum_{v\in
\bbZ_n}f_v(x)=0$ (these claims are also  easily verified directly).
Part (2) is obvious from the observation that $e_{v+1}-e_v=(x-1)x^v$, 
and to see (3), simply note that from (1), after
multiplication by $x-1$, we obtain $g_0(x)+\cdots+g_{n-1}=0$.  \epf

We can now derive an expression for the Smith group~$\Gamma(\gD)$ in
terms of polynomials.  
\ble{LLSmith}
For the Smith group $\Gamma(\gD)$ we have
\[ \Gamma(\gD)=\bbZ^n/\bbZ^n\gD=(\bbZ[x] \bmod x^n-1)/\langle f_v(x)
\mid 0\leq v\leq n-1\rangle_{\bbZ} = \bbZ\oplus \cZ_n /\langle f_v(x) \mid
1\leq v\leq n-1\rangle_{\bbZ}.\] 
\ele 
\bpf First note that
the vectors $\bbZ^n$ correspond to the polynomials in $\bbZ[x]\bmod
x^n-1$ and, since the rows of $\gD$ correspond to the polynomials
$f_v(x)$ or $-f_v(x)$, the vectors in the row space $\bbZ^n\gD$
correspond to the elements of $\langle f_v(x) \mid 0\leq v\leq n-1\rangle_{\bbZ}$. This shows the
second equality in the lemma.  Then, by 
Lemma~\ref{LLbasic} (2) and the Chinese Remainder Theorem, we have
that $\bbZ_n[x]\bmod x^n-1 \cong \bbZ\oplus \cZ_n$. Again by
Lemma~\ref{LLbasic}, every $f_v(x)$ is contained in~$\cZ_n$ and
$f_0(x)$ depends on the other $f_v(x)$, so the lemma follows.  \epf

The next result is one of the key points in our approach.
\btm{LTspan}
Let $n$ and $d$ be integers, with $n\geq 1$. If $|d|\geq
2$, then the polynomials $\gre_v$, for $1\leq v\leq n-1$  
(resp., the polynomials $g_v(x)$ ($1\leq v\leq n-1$)) are independent
over $\bbQ$,  and $\dim_{\bbQ} \cE_{n,d}=n-1$.
\etm 
\bpf
In view of Lemma~\ref{LLbasic} (3), it suffices to
show that the $g_v(x)$ for $1\leq v\leq n-1$ are
independent over~$\bbQ$. To see this, suppose that 
\[0\bmod x^n-1=\sum_{v\neq 0} a_vg_v(x)=d(x-1)\sum_{v\neq0} a_vx^v-
(x^d-1)\sum_{v\neq0} a_vx^{dv},\quad a_v\in\bbQ.\]
Writing $a(x)=\sum_{v>0} a_v x^v$ and
$c(x)=\sum_{v> 0} a_vx^v(x-1)=(x-1)a(x)
=\sum_{i} c_ix^i$,
we have that $c(x^d)=dc(x)\bmod x^n-1$. Subsituting $x$ with $x^d$, one obtains
$c((x^d)^d)=c(x^{d^2})=dc(x^d)=d^2c(x)$. Similarly, 
$c(x^{d^3})=dc(x^{d^2})=d^3c(x)$, etc.  Hence $c(x^{d^e})=d^ec(x)$  in
$\bbQ[x]\bmod x^n-1$, for every integer $e\geq 1$. However, if
$|d|\geq 2$ and $c(x)\neq 0$, then the left-hand side has bounded
coefficients while the right-hand side has an unbounded coefficient
$d^ec_i$ if $c_i\neq 0$, a contradiction. It follows that
$c(x)=(x-1)a(x)=0$, hence $a(x)\equiv 0\bmod 1+x+\cdots+x^{n-1}$.
But as $a(0)=0$ and $\deg a<n$, it follows
that $a_v=0$ for all $v$.  \epf
%
\bco{LCspan} For integers $d$ with $|d|\geq2$, the sandpile group (or
equivalently,  the critical group) $S(n,d)$ of the generalized de
Bruijn graph  $\mathrm{DB}(n,d)$ (if $d>0$) or of the generalized
Kautz graph $\mathrm{DB}(n,|d|)$ (if $d<0$) can be expressed as
\beql{LESPDB} 
S(n,d)=\cZ_n/ \langle f^{(n,d)}_v(x) \mid 1\leq v\leq n-1\rangle_{\bbZ}
\eeql
\eco 
\bpf First, we claim that the polynomials $f_v$ for $1\leq v\leq n-1$ 
are independent over~$\bbQ$. Indeed,
every nontrivial relation between the $f_v(x)$ implies (after
multiplication by $x-1$) a similar relation between the $g_v(x)$;
however, according to Theorem~\ref{LTspan}, such relation cannot exist
if $|d|\geq2$.  As a consequence, for 
$K:= \langle f_v(x) \mid 1\leq v\leq n-1 \rangle_{\bbZ}$ one has
 $\dim_\bbQ K =n-1$. 
In view of  what was stated in Section~\ref{ssect:scs}, this implies
that the quotient $\cZ_n/K$ in  Lemma~\ref{LLSmith}  is a
finite group, and the lemma follows.  \epf
It is not so easy to determine the structure of $S(n,d)$ by employing
(\ref{LESPDB}), due to the complicated form of the
polynomials~$f_v(x)$. The polynomials $g_v(x)=(x-1)f_v(x)$ have a much
easier structure, which motivates the following approach. We define
the {\em sand dune group\/} $\gS(n,d)$  of the generalized de Bruijn
graph $\mathrm{DB}(n,d)$  for $d\geq 2$ (resp. of the generalized Kautz
graph $\mathrm{Ktz}(n,|d|)$  for $d\leq -2$) as 
\[
\gS(n,d)=\cZ_n/\langle g_v(x)\mid 1\leq v\leq n-1\rangle_{\bbZ} =\cZ_n/\cE_{n,d}.
\]
%
The next result is crucial to our approach: it shows that 
$S(n,d)<\gS(n,d)$, and identifies the elements of 
$\gS(n,d)$ that are contained  in $S(n,d)$.
\btm{Lsub} The sand dune group $\gS(n,d)$ is finite, and the sandpile
group $S(n,d)$ is a subgroup of~$\gS(n,d)$.  Moreover, if
$a=\sum_{v=1}^{n-1} a_v e_v\in \gS(n,d)$, then $a\in S(n,d)$ 
if and only if  $\sum_{v=1}^{n-1} va_v\equiv0 \bmod n$.  
\etm 
\bpf 
The finiteness of $\gS(n,d)$ follows from the
fact that $\dim_\bbQ \cE_{n,d}=n-1$, as proved in Lemma~\ref{LTspan}.
Next,  write $T_v(x)=e_v(x)/(x-1)=1+\cdots+x^{v-1}$ for 
$0\leq v\leq n-1$.
Consider the map $\phi$ on $\bbQ[x]\bmod x^n-1$ 
for which $\phi(c(x))=(x-1)c(x)$. It is $\bbQ$-linear, and
$\ker\phi=\langle T_n(x)\rangle_{\bbQ}$. Since $T_n(1)=n\neq 0$, the
restriction of~$\phi$ to $\cZ_n$ is one-to-one, and $\phi$ maps
$\langle f_v(x) \mid 1\leq v\leq n-1\rangle_{\bbZ}$ onto
$\cE_{n,d}$. 
As $\phi(\cZ_n)\subseteq \cZ_n$, one sees that
$\phi$ embeds $S(n,d)$ as a subgroup
$\phi(\cZ_n)/\cE_{n,d}$ in $\gS(n,d)$. 
To determine that subgroup, we need to determine $\phi(\cZ_n)$.


To this end, let $a(x)=\sum_v a_ve_v(x) \in \cZ_n$. 
Then we can write
$a(x)=h(x)(x-1)$ with $h(x)=a(x)/(x-1)=\sum_v a_vT_v(x)$. Note that
since $T_v(1)=v$, we have $h(1)=\sum_v va_v$.
Now $a(x)=\phi(b(x))=b(x)(x-1)$ with $\deg b<n$ 
precisely when $b(x)$ is of the form $b(x)=h(x)-\gl T_n(x)$ with
$\gl\in \bbQ$; note that $b(x)\in \cZ_n$ precisely when $\gl\in \bbZ$
and $b(1)=h(1)-\gl T_n(1)=\sum_v va_v-\gl n=0$; such a $\gl$ exists
precisely when $\sum_v va_v\equiv 0\bmod n$.  
\epf 

\bco{LCord} 
We have $\gS(n,d)/S(n,d)=\bbZ_n$ and in particular
$|\gS(n,d)|=n|S(n,d)|$.
\eco 
\bpf The map $\theta:\bbQ[x]/(x^n-1)\to \bbZ_n$ given by
$\theta(\sum_v a_v e_v)=\sum_vva_v$ has the property that 
$\theta(\gre_v)=\theta(de_v-e_{dv})=0$. Hence
it is well-defined as a map on~$\gS(n,d)$;
it is obviously a homomorphism, and it is surjective  since
$\theta(e_v)=v$ for all $v\in \bbZ_n$.  As a consequence,
$n=|\bbZ_n|=|\im(\theta)|=|\gS(n,d)|/|\ker(\theta)|=|\gS(n,d)|/|S(n,d)|$.
\epf

We remark that the determination of $S(n,d)$ is complicated by the
fact that $S(n,d)$ is not always a {\em direct summand\/} of
$\gS(n,d)$,  as is illustrated by the following.  
\bex{LExnods}
Let $n=4$ and $d=3$. Then $\gS(n,d)=\bbZ_8\oplus \bbZ_2$ and
$S(n,d)=\bbZ_4$, which is not a direct summand of $\bbZ_8\oplus
\bbZ_2$.  \eex 
The above descriptions of $S(n,d)$ and 
$\gS(n,d)$, and the embedding of $S(n,d)<\gS(n,d)$, 
are quite suitable for the determination of these groups.
In that process, at several places  information is required about the
order of various group elements. Our next few results provide that
information.  
%
\ble{Lord} 
Every element $\ga\in \gS(n,d)$ can be expressed as
$\ga=\sum_{v>0} \ga_v\gre_v$, with $\ga_v\in\bbQ$ 
satisfying $0\leq a_v<1$ for each $1\leq v\leq n-1$; then the order of $\ga$ in
$\gS(n,d)$ is the smallest positive integer~$m$ such that $m\ga_v\in
\bbZ$ for each $1\leq v\leq n-1$.  
\ele 
\bpf 
According to Theorem~\ref{LTspan}, the $\gre_v$ are independent in $\bbQ[x]\bmod
x^n-1$. Therefore,  every polynomial $f(x)$ in $\bbQ[x]$
with $f(1)=0$ and $\deg f<n$ has a unique expression $f=\sum_{v>0}f_v\gre_v$ as
linear combination of the $\gre_v$. Such an expression is 0 modulo
$\cE_{n,d}$ if and only if  all coefficients $f_v$ are integers. Now the
claim is obvious.  
\epf

In order to use this result, we must be able to express the
polynomials $e_v$ in terms of the~$\gre_v$. This can be done as
follows.
\bde{LDfe} 
Let $v\in \bbZ_n$. Given $d\in\bbZ$, there are unique 
$\bbZ\ni e>0$ and $\bbZ\ni f\geq0$ such that the $d^iv$ in~$\bbZ_n$ are
distinct for $0\leq i\leq e+f-1$, while $d^{e+f}v=d^fv$. 
We say that $v$ has $d$-type $[f,e]$ in~$\bbZ_n$.
\ede
%
\ble{Lexp} Let $n$ and  $d$ be integers with $n\geq1$ and $|d|\geq2$.
If $v$ has $d$-type $[f,e]$ in $\bbZ_n$ then in~$\bbQ[x]\bmod x^n-1$,
we have \[e_v =\sum_{i=0}^{f-1}
\frac{1}{d^{i+1}}\gre_{d^iv}+\sum_{j=0}^{e-1}\frac{d^j}{d^f(d^e-1)}\gre_{d^{f+j}v}.\]
\ele 
\bpf 
First note that by a simple ``telescoping'' summation
\[e_v = d^{-1}\gre_{v}+\cdots+
d^{-f}\gre_{d^{f-1}v}+d^{-f}e_{d^fv}.\] 
Put $w=d^fv$.  Then 
\[e_w = d^{-1}\gre_{w}+\cdots+ d^{-e}\gre_{d^{e-1}w}+d^{-e}e_{d^ew}, \] 
and since $d^ew=w$, we conclude that \[(d^e-1)e_w = d^{e-1}\gre_w+\cdots +
\gre_{d^{e-1}w}.\] By combining these two results, the lemma follows.
\epf
%
%
In view of Lemma~\ref{Lord}, we immediately obtain the following.
\bco{Cord} Let $|d|\geq2$. If $v$ has $d$-type $[f,e]$ in $\bbZ_n$
then $e_v$ has order $|d^f(d^e-1)|$ in $\gS(n,d)$.  \eco 
%

\bre{LRjoint} \rm As we have seen, the group $S(n,d)$ is equal to the
sandpile group of~$\mathrm{DB}(n,d)$ if $d\geq2$  (resp.
of~$\mathrm{Ktz}(n,|d|)$ if $d\leq-2$).  In what
follows, we concentrate on the generalized de Bruijn graph and
therefore we assume $d\geq2$. We leave it to the reader to make
the necessary adaptations for the generalized Kautz graphs.
\ere

\section{\label{sect:main}Main results}
%
Let $n$ and $d$ be fixed integers with $n\geq1$ and $|d|\geq 2$ . The
description of the sandpile group
$S(n,d)$ and the sand-dune group $\gS(n,d)$ 
involves a sequence of numbers defined as follows.  Put
$n_0=n$, and for $i=1, 2,\ldots$, define $d_i= (n_i,|d|)$ and
$n_{i+1}=n_i/d_i$.  We have $n_0>\cdots>n_k=n_{k+1}$, where $k\geq0$
is the smallest integer for which $d_k=1$.  We refer to the
sequence $n_0>\cdots >n_k=n_{k+1}$  as the {\em $d$-sequence\/} of
$n$.  In what follows, we write $m=n_k$.
Note that $n=d_0\cdots d_{k-1}  m$ with $(m,d)=1$.

Since $(m,d)=1$, the map $x\mapsto dx$ is invertible and 
partitions $\bbZ_m$ into orbits  of the form $O(v)=\{v,dv,\ldots, d^{o(v)-1}v\}$. 
Here, $O(v)$ is sometimes referred to as the {\em $d$-ary cyclotomic coset
of~$v$ modulo $m$\/}.  We refer to $o(v)=|O(v)|$ as the 
{\em order\/} of $v$. 

For every prime $p|m$, we define $\pi_p(m)$ to be the largest power of
$p$ dividing $m$. Let $\cV$ be a set of representatives of
the orbits $O(v)$ different from $\{0\}$, where we ensure that for
every prime divisor $p$ of $m$, all integer numbers 
of the form $m/p^j$ are contained in~$\cV$.
(This is possible since no two of these numbers are in the same
cyclotomic coset.) 

\btm{LTmain} 
Let $n=d_0\cdots d_{k-1}  m$ with $(m,d)=1$.
The groups $S(n,d)$ and $\gS(n,d)$
are the sandpile and sand dune group of the generalized de Bruijn
graph $\mathrm{DB}(n,d)$ if $d\geq2$ (resp. 
of the generalized Kautz graph $\mathrm{Ktz}(n,|d|)$ if
$d\leq-2$).  With the above definitions and notation,
\beql{LEdb-sd} \gS(n,d)= \biggl[ \bigoplus_{i=0}^{k-1}
\bbZ_{|d|^{i+1}}^{n_i-2n_{i+1}+n_{i+2}}\biggl] \oplus \biggl[
\bigoplus_{v\in \cV}\bbZ_{|d^{o(v)}-1|}  \biggr] 
\eeql 
and
\beql{LEdb-sp}S(n,d)= \biggl[ \bigoplus_{i=0}^{k-1}
\bbZ_{|d|^{i+1}/d_i}\oplus
\bbZ_{|d|^{i+1}}^{n_i-2n_{i+1}+n_{i+2}-1}\biggl] \oplus \left[
\bigoplus_{v\in \cV} \bbZ_{|d^{o(v)}-1|/c(v)}\right], 
\eeql 
where, for each prime $p\mid m$, 
\begin{equation*}
c(v)=
\begin{cases}
\pi_p(m) & v=m/\pi_p(m), \text{ $p\neq 2$ or $d \equiv 1 \bmod
4$ or $4\nmid m$ }\\ 
\pi_2(m)/2 & v=m/\pi_2(m),  \text{ $d\equiv 3 \bmod 4$, and
$4\mid m$}\\
2 & v=m/2, \text{ $4\mid m$ and $d\equiv 3\bmod 4$}\\
1 & \text{ otherwise}.
\end{cases}
\end{equation*}
\etm 
Remark that since $n=d_0\cdots d_{k-1} m$
with $m=\prod_{p|m} \pi_p(m)$, the above result implies that 
$\gS(n,d)/S(n,d)=\bbZ_n$, in accordance with  the results in
Section~\ref{sect:prelims}.

%

With the notation from Section~\ref{ssect:Icm}, we have the following
isomorphisms, connecting critical groups and circulant matrices.
\btm{LTmain-circ} 
Let $d>0$ be a prime. Then $\gS(n,d)\cong C'(n,d)$ and 
$S(n,d)\cong C'(n,d)/\langle X\rangle$.  For $d$ a proper prime power, this
result also holds if $(n,d)=1$, but not always if $(n,d)\neq 1$.
\etm

The above results are proved in a number of steps. In what
follows, we outline the method for the generalized de Bruijn graphs;
for the generalized Kautz graphs, a similar approach can be used. 
First, we investigate  the ``multiplication-by-$d$'' map $d$
given by
$x\mapsto dx$ on the sandpile and sand dune groups. Let $\gS_0(n,d)$
and $S_0(n,d)$ denote the kernel of the map $x\mapsto d^kx$ on $\gS(n,d)$ and
$S(n,d)$, respectively.  It is not difficult to see that
$\gS(n,d)\cong \gS_0(n,d)\oplus \gS(m,d)$ and $S(n,d)\cong S_0(n,d)
\oplus S(m,d)$. Then, we use the map $d$ to determine $\gS_0(n,d)$ and
$S_0(n,d)$. It is easy to see that for {\em any\/} $n$, we have
$d\gS(n,d)\cong \gS(n/(n,d), d)$ and $dS(n,d)\cong S(n/(n,d), d)$.
With some more effort, it can be shows that the kernel of the map
$d$
on $\gS(n,d)$ (resp. on $S(n,d)$) is isomorphic to $\bbZ_d^{n-n/(n,d)}$ 
(resp. to $\bbZ_{d/(n,d)}\oplus \bbZ_d^{n-1-n/(n,d)}$). Then we use
induction on the length $k+1$ of the $d$-sequence of~$n$ to show
that $\gS_0(n,d)$ and $S_0(n,d)$ have the form of the left hand parts
of (\ref{LEdb-sd}) and (\ref{LEdb-sp}), respectively. This part of the
proof, although much more complicated, resembles the method used by
\cite{LeLI-sand} and \cite{BiKi-knuth}. 

Then it remains to handle the parts $\gS(m,d)$ and $S(m,d)$, where
$(m,d)=1$.  For the ``helper'' group $\gS(m,d)$ that embeds $S(m,d)$,
this is trivial: it is easily seen that $\gS(m,d)= \oplus_{v\in \cV}
\langle e_v\rangle$, and the order of $e_v$ is equal to the size
$o(v)$ of its orbit $O(v)$ under the map $d$, so  (\ref{LEdb-sd})
follows immediately. The $e_v$ are not contained in~$S(m,d)$, but we
can try to modify them slightly to obtain a similar decomposition for
$S(m,d)$. The idea is to replace $e_v$ by a modified version
$\te_v=e_v-\sum_{p|m} \gl_p(v) e_{m \pi_p(v)/\pi_p(m)}$, 
where the numbers $\gl_p(v)$ are chosen such that $\te_v\in S(m,d)$,
or by a suitable multiple of $e_v$, in some exceptional cases (these
are the cases where $c(v)>1$). It turns out that this is indeed
possible, and in this way the proof of Theorem~\ref{LTmain} can be
completed.

The proof of Theorem~\ref{LTmain-circ} is by reducing to the case
$(n,p)=1$ by an explicit construction, and then by diagonalizing
$C(n,p)$ over an appropriate extension of $\mathbb{F}_p$. Essentially,
as soon as $(n,p)=1$, one can read off a decomposition of $C(n,p)$
into cyclic factors from the irreducible factors of the polynomial
$x^n-1$ over $\mathbb{F}_p$. 

In the next sections, we provide the details of the proofs as
outlined above.

\section{The multiplication-by-$d$ map}

In the remainder of this section, we use the
map $x\mapsto dx$ on~$\gS(n,d)$ to determine
the structure of $\gS_0(n,d)$ and $S_0(n,d)$, i.e. the kernels of the map
$x\mapsto d^kx$. We require the following
simple result.

\btm{LTd} 
For any pair $n,d$ of positive integers,
$d\gS(n,d) \cong \gS(n/(n,d),d)$ and $dS(n,d)=S(n/(n,d),d)$.  
\etm
\bpf
Write $n_1=n/(n,d)$.
Define $\varphi:d\gS(n,d)\to \gS(n_1,d)$ by 
$\varphi(de_v)=e_{v \bmod n_1}$ for $1\leq v\leq n-1$
and extend $\varphi$ by linearity. We claim
that $\varphi$ defines an isomorphism between $d\gS(n,d)$ and
$\gS(n_1,d)$. To see this, proceed as follows. 
Since $de_v=e_{dv}$ in $\gS(n,d)$, we have that  
$d \sum_{v>0} \ga_v e_v= \sum_{v>0} \ga_v e_{dv}$, and from
Lemma~\ref{Lexp} we conclude that $\sum_{v>0} \ga_v e_{dv}$  can be
expressed as a linear combination of elements $\gre_{dv}$
in~$\cE_{n,d}$ with rational coefficients; the expression is 0
in~$\gS(n,d)$ if and only if all coefficients can be chosen to be
integer. So, noting that $\varphi$ maps $\gre_{dv}$ in $\cE_{n,d}$ to
$\gre_{v \bmod n_1}$ in~$\cE_{n_1,d}$, we conclude that $\varphi$ is
well-defined and in fact one-to-one on~$d\gS(n,d)$. Since $\varphi$ is
obviously onto $\gS(n_1,d)$, the desired conclusion follows.

To see that $\varphi$ also induces an isomorphism between $dS(n,d)$
and $S(n_1,d)$, in view of Theorem~\ref{Lsub} it is sufficient to
remark that for an element 
$d\sum_{v>0} a_ve_v=\sum_{v>0} a_ve_{dv}\in d\gS(n,d)$, 
we have $\sum_{v>0} dv a_v\equiv 0
\bmod n$ if and only if $\sum_{v>0} va_v\equiv 0\bmod n_1$.
\epf

The next step is to determine the {\em kernel\/} of the
multiplication-by-$d$  map $d:\bbQ[x]/(x^n-1)\to \bbQ[x]/(x^n-1)$,
defined by $x\mapsto dx$,
on both $\gS(n,d)$ and $S(n,d)$.  The
result is as follows.  
\btm{LTker} 
(i) The kernel $\ker_\gS(d)$ of the
map $d$ on $\gS(n,d)$ is isomorphic to
$\bbZ_d^{n-n_1}$.\\ 
(ii) The  kernel $\ker_S(d)$ of the map $d$ 
on $S(n,d)$ is isomorphic to $\bbZ_{d/d_0}\oplus \bbZ_d^{n-n_1}$.
\etm 
\bpf 
The order of $\gS(n,d)$ is equal to the
product of its invariant factors, which are the positive invariant
factors of the  $(n-1)\times (n-1)$ matrix $\gS=\gS^{(n,d)}$ that has
as rows the vectors 
$\gre_v=de_v-e_{dv}$ with respect to the basis 
$\{e_v \mid 1\leq v\leq n-1\}$.
Since $\Sigma$ is nonsingular, this product is equal to 
$\det\Sigma$. Now partition the 
set $\bbZ_n\setminus \{0\}$ of row and column indices of $\Sigma$
into parts
$\bbZ-d\bbZ_n$ and $(d\bbZ_n)\setminus \{0\}$. Under this
ordering of the rows and columns,  $\Sigma$ takes the form
\[
\gS=\gS^{(n,d)}= \left(\begin{array}{c|c} D & A\\ \hline
O & \gS', \end{array} \right),
\] 
where $D=\diag(d, \ldots, d)$ is a $(n-n_1)\times (n-n_1)$ diagonal matrix 
and $\gS'=\gS^{(n_1,d)}$ is the matrix corresponding to
$\gS(n_1,d)$. (Note that $d\bbZ_n\cong \bbZ_{n/(d,n)}=\bbZ_{n_1}$.)
We conclude that 
\beql{LEord} |\gS(n,d)| = d^{n-n_1}|\gS(n_1,d)|.
\eeql
In view of Theorem~\ref{LTd}, we have that
$|\ker_\gS(d)|=|\gS(n,d)|/|\im_\gS(d)|=d^{n-n_1}$, and as a
consequence of Corollary~\ref{LCord}, we also have 
\[|\ker_S(d)|
= |S(n,d)|/|S(n_1,d)|=(|\gS(n,d)|/n)/(|\gS(n,d)|/n_1)=d^{n-n_1}/d_0.\]
To actually construct a basis for these kernels, let us define 
\[\gD_{ab}:=e_{a+bn_1}-e_a = d^{-1}(\gre_{a+bn_1}-\gre_a),\] 
where the second equality follows directly from \eqref{def:e_v}. 
By Lemma~\ref{Lord}, each $\gD_{ab}$ has order~$d$. Hence they are 
contained in $\ker_\gS(d)$.  First, we claim that the
set 
\[ \cB=\{\gD_{ab} \mid 0\leq a\leq n_1-1, 1\leq b\leq d_0-1\}\] 
is independent in~$\gS(n,d)$ and a basis for $\ker_\gS(d)$.
Indeed, consider a $\bbZ$-linear combination of the elements of~$\cB$.  Since
\[
\sum_{a=0}^{n_1-1}\sum_{b=1}^{d_0-1}\gl_{ab}\gD_{ab} =
\sum_{a=0}^{n_1-1}\sum_{b=1}^{d_0-1} \gl_{ab}d^{-1}(\gre_{a+bn_1}
-\gre_a) = d^{-1}\sum_{a=0}^{n_1-1}\left(
\sum_{b=1}^{d_0-1}\gl_{a,b}\gre_{a+bn_1}-\biggl(\sum_{b=1}^{d_0-1}
\gl_{ab}\biggr)\gre_a \right), 
\] 
and since $a+bn_1=a'+b'n_1$ with
$a,a'\in\{0,1,\ldots, n_1-1\}$ and $b,b'\in\{0, \ldots, d_0-1\}$ is
only possible when $(a,b)=(a',b')$, each $\gre_{a+bn_1}$ occurs only
once in the expression. Hence the linear combination can be zero only
if every term $d^{-1}\gl_{ab}\gre_{a+bn_1}$ is zero, i.e. 
$d\mid \gl_{ab}$, i.e.  $\gl_{ab}\gD_{ab}=0$.

Since every $\gD_{ab}$ has order $d$, we conclude that
\[\ker_\gS(d)=\bigoplus_{\stackrel{a\in\{0,\ldots, n_1-1\}}{b\in\{1,
\ldots, d_0-1\}}} \langle \gD_{ab}\rangle \cong \bbZ^{n-n_1},\] 
where
the equality (instead of a containment) follows from the equality of
the respective sizes.

Similarly, the set 
\[ B=\{\gD_{ab}-b\gD_{01} \mid 0\leq a\leq n_1-1, 
1\leq b\leq  d_0-1, (a,b)\neq (0,1)\} \cup
\{d_0\gD_{01}\}\] spans $\ker_S(d)$: according to Theorem~\ref{Lsub},
every element of~$B$ is contained in~$S(n,d)$, and their independence
easily follows from the independence of $\cB$; a counting argument
similar to the one above shows that they span the entire kernel.
\epf

In what follows we need a few simple properties of finite
abelian groups. The first of these is a straightforward consequence of
the uniqueness of decomposition of finite abelian groups into cyclic groups of
prime power order.  
\bpn{LTdiv}
Let $G, H, K$ be finite abelian groups. 
If $G\oplus K\cong H\oplus K$, then $G\cong H$.
\qed
\epn 
The next result are needed when we deal with invariant
factors of a group. Recall that as a consequence of the uniqueness of
the Smith Normal Form, every abelian group $G$ also has a unique
decomposition $G\cong \bbZ_{s_1}\oplus \cdots \oplus\bbZ_{s_r}$ with
$1<s_1|\cdots |s_r$.  We refer to $s_1, \ldots, s_r$ as the {\em
invariant factors\/} of~$G$.  \btm{LTinv} (i) We have that
$\bbZ_m\oplus \bbZ_n\cong \bbZ_{(m,n)}\oplus \bbZ_{[m,n]}$, where
$(m,n)$ and $[m,n]$ denote the greatest common divisor (gcd) and the
least common multiple (lcm) of $m$ and $n$.\\ (ii) 
If $G$ has invariant factors $s_1, \ldots, s_r$, then $G\oplus \bbZ_m$
has invariant factors $s_1', \ldots, s_{r+1}'$ with $s_1'=(s_1, m)$,
$s_i'=(s_i, [s_{i-1},m])$ for $2 \leq i\leq r$, and
$s_{r+1}'=[s_r,m]$, or invariant factors $s_2', \ldots, s_{r+1}'$ if
$s_1'=(s_1,m)=1$.  \etm \bpf Part (i) is trivial, it follows
immediately from the decomposition of $\bbZ_m$ and $\bbZ_n$ into their
prime power summands. Part (ii)  follows from part (i)  and induction
on~$r$: we have that 
\beq  
\bbZ_{ s_1}\oplus \cdots \oplus \bbZ_{s_r} \oplus \bbZ_m \cong  \bbZ_{
s_1'}\oplus \cdots \oplus \bbZ_{s_{r-1}'}\oplus\bbZ_{[s_{r-1},m]}
\oplus \bbZ_{s_r} \cong  \bbZ_{ s_1'}\oplus \cdots \oplus
\bbZ_{s_{r-1}'}\oplus\bbZ_{s_{r}'} \oplus \bbZ_{s_{r+1}'}, 
\eeq 
where we have used part (i) and the fact that $s_{r-1}|s_r$. Since
$s_1'|\cdots |s_{r+1}'$ and $1<s_1|s_2'$, the invariant factors of
$G\oplus \bbZ_m$ are $s_1', \ldots, s_{r+1}'$, or  $s_2', \ldots,
s_{r+1}'$ in case that $s_1'=(s_1,m)=1$.  \epf
We also require the following simple lemma.  
\ble{LLpre} Let $G$ be an
abelian group, and let $d, m$ be two positive integers. Let $d'$
denote the maximal divisor of $d$ for which $(d',m)=1$. If $\bbZ_m$ is
a direct summand of $dG$, then $\bbZ_{md/d'}$ is a direct summand
of~$G$.
In particular, if $d|m$, then $\bbZ_{md}$ is a direct summand of~$G$.
\ele 
\bpf 
Consider a decomposition of $G$ into cyclic
groups of prime-power order. If $\bbZ_{p^t}$ is a direct summand
of~$G$ and if $p^s||d$, then
$d\bbZ_{p^t}=\bbZ_{p^t/(d,p^t)}=\bbZ_{p^{max(0,t-s)}}$ is the
corresponding direct summand in~$dG$. Therefore, if $p^r||m$, then
the direct summand $\bbZ_{p^r}$ of~$\bbZ_m\leq dG$ can only arise
from a direct summand $\bbZ_{p^{r+s}}$ of $G$, i.e. the required
direct summand of $\bbZ_{md/d'}$.
\epf
%

Let $\gd=d_0\cdots d_{k-1}$
and write $m=n_k$. Then $n=\gd m$ with $(\gd,m)=1$.  The Chinese
Remainder Theorem (CRT)  decomposition $\bbZ_n=\bbZ_\gd \oplus
\bbZ_{m}$ induces a corresponding decompositions for the sand dune and
sandpile groups.  
\ble{LTcrt}We have that \[\gS(n,d)\cong
\gS_0(n,d)\oplus \gS(m,d)\] and \[S(n,d)\cong S_0(n,d) \oplus
S(m,d),\] where $\gS_0(n,d)=\ker_{\gS}(d^k)$  and $S_0(n,d)
=\ker_{S}(d^k)$ are the kernel of the map $x\mapsto d^kx$ on
$\gS(n,d)$ and $S(n,d)$, respectively.  
\ele 
\bpf Since $n=\gd m$ with
$(\gd,m)=1$, by the CRT there are integers $\chi$ and $\eta$ such that
$\chi \gd$ and $\mu m$ are mutually orthogonal idempotents, that is,
$\chi \gd \equiv  1 \bmod m$ and $\mu m\equiv 1 \bmod \gd$. As a
consequence, for each $v\in \bbZ_n$ we have $v=(v\mu m) + (v\chi
\gd)$, and it is easily seen that the map $v\mapsto (v\mu m, v\chi
\gd)$ induces a decomposition $\bbZ_n \cong \bbZ_\gd  \oplus
\bbZ_{m}$. Then we can write \[e_v = (e_{v\chi \gd+v\mu m} - e_{v\chi
\gd}) + e_{v\chi \gd},\] and it is easily verified that the map \[ e_v
\mapsto (e_{v\chi \gd +v\mu m} - e_{v\chi \gd}, e_{v\chi \gd})\]
induces a decomposition \[\gS(n,d)\cong \gS_0(n,d)\oplus \gS(m,d),\]
where $\gS_0(n,d)$ denotes the subgroup generated by the elements
$e_{v\chi \gd+v\mu m} - e_{v\chi \gd}$ for $v\in \bbZ_n$.  Now
$(d,m)=1$, so the map $x\mapsto dx$ acts as a permutation
on~$\bbZ_m$; since $de_v=e_{dv}$ on $\gS(n,d)$, we conclude that
$d\gS(m,d)\cong\gS(m,d)$. Next, since \[\gd e_{v\chi \gd + v\mu m} =
e_{v\chi \gd^2 + v\mu m\gd}= e_{v\chi \gd^2}=\gd e_{v\chi \gd},\] we
have that $\gd (e_{v\chi \gd + v\mu m} - e_{v\chi \gd})=0$. Now
$d_i=(d,n_i)|d$ for $i=0, \ldots, k-1$, so that $\gd|d^k$. Combining
these observations, we conclude that
$\gS_0(n,d)=\ker_{\gs(n,d)}(d^k)$.

By Theorem~\ref{Lsub},  the element 
$a=\sum_{v\geq 0}^{n-1}a_ve_v\in S(n,d)$ iff 
$\chi:=\sum_{v\geq 0}{n-1}va_v\equiv0 \bmod n$, which by CRT is the case if
and only if  $\chi$ is 0 both modulo $\gd$ and
modulo $m$. Therefore, $a\in S(n,d)$ iff both projections $\sum_v
a_v (e_{v\chi \gd +v\mu m} - e_{v\chi g})$  and $\sum_v a_v e_{v\chi
\gd}$ are in $S(n,d)$. It follows that the above decomposition for
$\gS(n,d)$ induces a similar decomposition for $S(n,d)$.  \epf

We now use the multiplication-by-$d$ map to inductively determine
the parts $\gS_0(n,d)$ and $S_0(n,d)$.  
\btm{LTmain1} Let $n$ have
$d$-sequence $n=n_0, n_1, \ldots, n_k=n_{k+1}$, with
$d_i=n_i/n_{i+1}=(n_i,d)$ for $i=0, \ldots, k$. Then 
\[\gS_0(n,d) \cong\bigoplus_{i=0}^{k-1} \bbZ_{d^{i+1}}^{n_i-2n_{i+1}+n_{i+2}}\] and
\[S_0(n,d) \cong \bigoplus_{i=0}^{k-1}\biggl(\bbZ_{d^{i+1}/d_i}\oplus
\bbZ_{d^{i+1}}^{n_i-2n_{i+1}+n_{i+2}-1}\biggr).\] 
\etm 
\bpf 
We use induction on~$k$.  If $k=0$, then there is nothing to prove.  Now,
suppose that $k\geq1$.  By Theorem~\ref{LTd}, we have
$d\gS_0(n,d)\cong \gS_0(n_1,d)$.  Since $n_1$ has $d$-sequence $(n_1,
n_2, \ldots, n_k)$,  by induction and by Lemma~\ref{LLpre}, we
conclude that $\gS_0(n,d)$ is of the form \[\gS_0(n,d) \cong \Lambda
\oplus \bigoplus_{i=1}^{k-1} \bbZ_{d^{i+1}}^{n_i-2n_{i+1}+n_{i+2}}\]
with $\Lambda\subseteq \ker_\gS(d)$. 
Hence $\ker_\gS(d)=\ker_{\gS_0}(d)=\Lambda\oplus \bbZ_d^{n_1-n_2}$ (to
see this, note that
$n_i-2n_{i+1}+n_{i+2}=(n_i-n_{i+1})-(n_{i+1}-n_{i+2})$ for all $i$ and
$n_k-n_{k+1}=0$).  So by Theorem~\ref{LTker} and
Proposition~\ref{LTdiv}, we have that $\Lambda\cong
\bbZ_d^{(n-n_1)-(n_1-n_2)}$, as was to be proved.

Similarly, again using Lemma~\ref{LLpre}, we can conclude that
\[S_0(n,d) \cong L \oplus \bbZ_{d^2}^{n_1-2n_2+n_3-1} \oplus
\bigoplus_{i=2}^{k-1}\biggl(\bbZ_{d^{i+1}/d_i}\oplus
\bbZ_{d^{i+1}}^{n_i-2n_{i+1}+n_{i+2}-1}\biggr),\] where \beql{LEdL}
dL=\bbZ_{d/d_1}.\eeql
From this expression, we read off that
$\ker_S(d)=\ker_{S_0}(d)=\ker_L(d)\oplus \bbZ_d^{n-n_1-1}$, hence by
Theorem~\ref{LTker} and Proposition~\ref{LTdiv}, we conclude that
\beql{LEkL} \ker_L(d) = \bbZ_{d/d_0}\oplus \bbZ_d^{n-2n_1+n_2}.\eeql
Suppose that $L$ has invariant-factor decomposition
\[L=\bbZ_{s_0}\oplus \cdots \oplus \bbZ_{s_r}\] with $s_0|\cdots |
s_r$. Then from (\ref{LEdL}) and (\ref{LEkL}), we conclude that
\beql{LEdLp} \bbZ_{s_0/(s_0,d)} \oplus \cdots \oplus
\bbZ_{s_r/(s_r,d)} =\bbZ_{d/d_1}\eeql and \beql{LEkLp}\bbZ_{(s_0,d)}
\oplus \cdots \oplus \bbZ_{(s_r,d)} =\bbZ_{d/d_0}\oplus
\bbZ_d^{n-2n_1+n_2}.\eeql Since $(s_i,d)|(s_{i+1},d)$ for $i=0,
\ldots, r-1$ and $(d/d_0)| d$, both direct sums in (\ref{LEkLp}) are
invariant-factor decompositions. Hence $(s_0,d)=d/d_0$ and $(s_i,d)=d$
for $i=1, \ldots, r$. So we can write $s_0=\tau d/d_0$ for some $\tau$
with $(\tau,d_0)=1$ and $s_i=\gs_i d$ for $i=1, \ldots, r$, where
$\gs_1|\cdots |\gs_r$ and $\tau|\gs_1d_0$, so that $\tau|\gs_1$.
Moreover, from (\ref{LEdLp}), we conclude that \[\bbZ_{\tau}\oplus
\bbZ_{\gs_1}\oplus \cdots \oplus \bbZ_{\gs_r} = \bbZ_{d/d_1}.\] Now,
using Theorem~\ref{LTinv}, we conclude that the left-hand side above
has invariant factors $(\gs_1,\tau)$, $(\gs_{i+1}, [\gs_i,\tau])$ for
$i=1, \ldots, r-1$, and $[\gs_r,\tau]$, while the right-hand side has
invariant factors $d/d_1$. Since the invariant factors are unique, we
conclude that  $(\gs_1,\tau)=1$ and $(\gs_{i+1}, [\gs_i,\tau])=1$ for
$i=1, \ldots, r-1$, while $[\gs_{r},\tau]=d/d_1$. Since
$\gs_i=(\gs_{i+1}, \gs_i)| (\gs_{i+1}, [\gs_i,\tau])$, we conclude
that $\gs_1=\ldots = \gs_{r-1}=1$ and $(\gs_r,\tau)=1$, while
$[\gs_r,\tau] =\gs_r\tau=d/d_1$. Since $\tau|\gs_1|\gs_r$,
we conclude that $\tau=1$ and $\gs_r=d/d_1$,
hence $L=\bbZ_{d/d_0}\oplus \bbZ_d^{n-2n_1+n_2-1}\oplus
\bbZ_{d^2/d_1}$, which is what we wanted to prove.  \epf


\subsection{Adaptations for the case of generalized Kautz graphs}
%
With minor adaptations, all the results in this section are also valid
when $d<0$, so for the sand dune group $\gS(n,d)$ and sandpile group
$S(n,d)$ of the generalized Kautz graph ${\rm Ktz}(n,|d|)$. Like
before, write $n_1=n/(n,|d|)$ and $d_0=(n,|d|)$. Using multiplication
by~$d$, we conclude in a similar way that $d\gS(n,d)\cong \gS(n_1,d)$
and $dS(n,d)\cong S(n_1,d)$,  and also that $\ker_{\gS}(d)\cong
\bbZ_d^{n-n_1}$ and $\ker_{S}(d)\cong \bbZ_{|d|/d_0}\oplus
\bbZ_{|d|}^{n-n_1-1}$. And finally, we can use these facts to
determine $\gS_0(n,d)$ and $S_0(n,d)$ in a similar way.


\section{\label{sect:rp}The sand dune and sandpile group in the
relatively prime case}
%
In this section, we determine the sand dune group $\gS(m,d)$ and the
sandpile group $S(m,d)$ for fixed positive integers $m$ and $d$ with
$(m,d)=1$.  In that case, the map $x\mapsto dx$ partitions
$\bbZ_m$ into orbits of the form \[O(v) = \{v,dv, \ldots,
d^{e-1}v\},\] where $d^ev\equiv v\bmod m$ and $d^iv\not\equiv d \bmod
m$ for $1\leq i<e$. 
Recall that we refer to $o(v)=|O(v)|$ as the {\em order\/} of~$v$, and
that $\cV$  denotes a complete set of representatives of the orbits
different  from~$\{0\}$, that is, $\cV$ contains precisely one element
from each 
orbit different from~$\{0\}$. Recall that for $p$ a prime, $\pi_p(v)$
denote the largest power of~$p$ dividing~$v$. For the remainder of
this section, we let 
\[ \bbP=\{ p \mid \mbox{$p$ is a prime and $p|m$}\},\] 
and write 
\[M_p=m/\pi_p(m), \qquad p\in \bbP.\]
For later use, we define 
\[ \cV^* = \{1\leq v\leq m \mid \mbox{$v\equiv p^i M_p
\bmod m$ for some $p\in \bbP$ and some  integer $i\geq 0$}\}.\] Since
$(d,m)=1$, it is easily seen that for $p,q\in\bbP$, if $p^iM_p \equiv
d^k q^jM_q \bmod m$ then $p=q$ and $i=j$. Thus the
elements in $\cV^*$ are in different orbits on $\bbZ_m$. 
(Another way to see this is to note that every divisor $k|m$ is
minimal in its orbit, which is contained in $d\bbZ_m$.)

\smallskip \noindent {\em For the remainder of this section, we 
assume that the set $\cV$ of orbit representatives contains all the
members of~$\cV^*$.\/} \smallskip

The determination of the sand dune group is easy.  
\btm{TSd} With the
above notation, we have that
$\gS(m,d)=\bigoplus_{v\in\cV}\bbZ_{d^{o(v)}-1}$.  
\etm 
\bpf 
By Lemma~\ref{Lexp}, the expression for $e_v$ in terms of
the~$\gre_w$ involves only $\gre_w$ with $w\in O(v)$. Hence the $e_v$ with
$v\in \cV$ are independent. Moreover, since $d^ie_v=e_{d^iv}$, the
subgroup $\langle e_v\rangle$ generated by $e_v$ contains every $e_w$
with $w\in O(v)$, so $\gS(m,d) =\oplus_{v\in \cV}\langle e_v\rangle$,
According to Corollary~\ref{Cord}, 
the order of $e_v$ is equal to $d^{o(v)}-1$, so $\langle
e_v\rangle\cong \bbZ_{d^{o(v)}-1}$, from which the theorem follows.
\epf

\bre{LRSd} \rm Aa alternative way to see the above result is to remark
that  the matrix $\gS=\gS^{(m,d)}$ is equivalent to a block-diagonal
matrix with $|\cV|$ blocks $\gS_v$ ($v\in \cV$), where the block
$\gS_v$ is the restriction of $\gS$ to the rows and columns indexed by
orbit $O(v)$; moreover, if within an orbit $O(v)$ we index in the
order  $v,dv,\ldots, d^{o(v)-1}v$, then $\gS_v$ is $o(v)\times o(v)$
and of the form \[ \gS_v=
\left( \begin{array}{ccccc} d&-1&0&\cdots&0\\ 0&d&-1&\cdots&0\\
&\ddots&\ddots &\ddots&  \\ 0&&0&d&-1\\ -1&0&   &0 &   d \end{array}
\right).  \] Now it is easy to see that $\gS_v$ has invariant factors
$(1, \ldots, 1, d^{o(v)}-1)$, for example, by successively adding $d$
times the first row to the second row, then $d$ times the current
second row to the third row, \ldots, and finally  $d$ times the then
current $(o(v)-1)$th row to the last row, and then subtracting from
the first column a suitable linear combination of the other columns.
\ere



As explained in Section~\ref{sect:main}, it is a bit more complicated
to determine the structure of~$S(m,d)$. First, we define the
modified generators $\te_v$ for $\gS(m,d)$. To this end, we need some
preparation.  Let $p\in \bbP$. Since $\pi_p(m)$ and $M_p=m/\pi_p(m)$
are relatively prime,  there is a number $\eta_p$ such that \[\eta_p
M_p  \equiv 1 \bmod \pi_p(m).\] Define \[\gl_p(v)=\eta_p v/\pi_p(v),\]
and for $v\in \cV\setminus \cV^*$, put \[\te_v=e_v-\sum_{p|m} \gl_p(v)
e_{\pi_p(v)M_p}.\]

\btm{LTte1} For $v\in \cV\setminus \cV^*$,  let $\te_v $ be defined as
above. Then $\te_v$  is contained in~$S(m,d)$, and $\te_v$ and $e_v$
have the same the order in~$S(m,d)$. Moreover, the $\te_v$ for $v\in
\cV\setminus \cV^*$ together with the $e_v$ for $v\in \cV^*$ are
independent and generate $\gS(m,d)$.  \etm \bpf To show that $\te_v\in
S(m,d)$, we use Theorem~\ref{Lsub}.  For a fixed prime
$q\in\bbP$, the numbers $M_p=m/\pi_p(m)$ with $p\in \bbP\setminus
\{q\}$ are divisible by $\pi_q(n)$, so that  modulo~$\pi_q(v)$, we
have \[v-\sum_{p|m}  \gl_p(v) \pi_p(v)M_p\equiv v-\gl_q(v)\pi_q(v)M_q
\equiv v-\eta_q (v/\pi_q(v))\pi_q(v)M_q\equiv 0 \bmod \pi_q(v).\]
Since this holds for every $q\in \bbP$, by the Chinese Remainder
Theorem we have that $v-\sum_{p|m} \pi_p(v) \pi_p(v)m/\pi_p \equiv 0
\bmod m$, hence by Theorem~\ref{Lsub}, $\te_v$ is in $S(m,d)$.

Next, recall that by Corollary~\ref{Cord}, $e_v$ has order
$d^{o(v)}-1$. To show that $\te_v$ has order $d^{o(v)}-1$, it is
sufficient 
to show that $(d^{o(v)}-1)e_{\pi_p(v)m/\pi_p(m)}=0$ holds for every
$p\in \bbP$. To see this, we proceed as follows. By the definition of
$o(v)$, we  have that  $(d^{o(v)}-1)v\equiv 0\bmod m$, hence
$(d^{o(v)}-1)\pi_p(v) \equiv 0 \bmod \pi_p(m)$, and therefore
$(d^{o(v)}-1)\pi_p(v) m/\pi_p(m) \equiv 0 \bmod m$.  So the order of
$e_{\pi_p(v)M_p}$ divides $d^{o(v)}-1$, 
and now the desired conclusion follows from the
equality $d^{o(v)}e_{\pi_p(v)M_p}= e_{d^{o(v)}\pi_p(v)M_p}$.

Finally, by the definition of $\cV^*$, it is obvious that the $\te_v$
for $v\in \cV\setminus \cV^*$ together with the $e_v$ for $v\in \cV^*$
have the same span as the $e_v$ for $v\in \cV$, from which the claim
follows immediately.  \epf

So now we are left with the choice of suitable elements $\te_v$ for
$v\in \cVs$.
First, we need a simple number-theoretic result. For a prime $p$, let
$\nu_p(n)$ denote the largest integer~$e\geq 0$ for which
$p^e|n$. We say that an integer $d$ has {\em order $e$ modulo $p$\/}
if $e$ is the smallest positive integer for which $p|d^e-1$.
\ble{Lpcont} For a prime $p$,with $(p,d)=1$, let $d$ have order $e$
modulo $p$, and suppose that $\nu_p(d^e-1)=a$.
Then $d$ has order $ep^i$ modulo $p^{a+i}$ for all $i\geq0$, except
when $p=2$ and $d\equiv 3 \bmod 4$. In that exceptional case, $e=1$
and $a=1$, and if $\nu_2(d^{2}-1)=b$, then $b\geq3$ and $d$ has order
$2^{i+1}$ modulo $2^{b+i}$ for all $i\geq0$ (and order $1$ modulo
$2^i$ for $1\leq i<b$).  \ele \bpf For an integer $t\geq1$, if $d$ has
order $f$ mod $p^t$, then $p^t|d^n-1$ if and only if $f|n$. So the
order of $d$ modulo $p^{t+1}$ is of the form $kf$.  Moreover, if
$\nu_p(d^f-1)=s\geq t$, then $d^f=1+qp^s$ for some integer $q$ not
divisible by~$p$, so \[ d^{rf}=(1+qp^s)^r\equiv 1+rqp^s \bmod
p^{2s};\] hence for every integer $k\geq1$,
\[d^{kf}-1=(d^f-1)(1+d^f+\cdots +d^{(k-1)f}) \equiv 
k+qp^a(0+1+\cdots +(k-1)) \equiv qp^s \left(k+qp^s\binom{k}{2}\right) \bmod
p^{3s}.\] 
So the smallest $k>1$ for which $p^{s+1}|d^{kf}-1$ is $k=p$,
in which case \[d^{pf}-1\equiv qp^{s+1}(1+qp^{s-1}\binom{p}{2}) \bmod
p^{3s}.\] Moreover, we see that $\nu_p(d^{pf}-1)=s+1$ , except when
$s=1$ and $p=2$. In that case, $q$ is odd and $8|d^{2e}-1$, and
$d\equiv 3 \bmod 4$ since $s=1$. So we conclude that there is a
``jump'' in the order of $d$ modulo powers $p^s$ of $p$ if and only if
$s=1$, $p=2$, and $d\equiv 3 \bmod 4$, as claimed in the theorem.
\epf

\bco{LLordcon} If $p\neq 2$ or $d\equiv 1\bmod 4$, then
$\nu_p(d^{o(M_p)}-1)-\nu_p(d^{o(p^tM_p)}-1) \leq t$ for $1\leq t \leq
\pi_p(m)-1$.  Otherwise, i.e. for $p=2$ and $d\equiv 3\bmod
4$, we have that $\nu_p(d^{o(M_p)}-1)-\nu_p(d^{o(p^tM_p)}-1) \leq t$
for $1\leq t \leq \pi_p(m)-2$.  
\eco 
\bpf Let $p\in \bbP$, and write
$s=\nu_p(m)$, so that $M_p=m/p^s$.  First we claim that the order of
$d$ modulo $p^{s-t}$ is equal to $o(p^tM_p)$.   To see this, note that
by definition,  $o(p^tM_p)$ is the smallest integer $e\geq 1$ for
which $(d^e-1)p^tM_p \equiv 0\bmod m$, or, equivalently, for which
$(d^e-1)p^t \equiv 0\bmod p^s$, from which the claim follows. 

Now suppose that $o(p^{s-1}M_p)=e$ and $\nu_p(d^e-1)=a$.  Then we see
from  Lemma~\ref{Lpcont} that in the ``non-exceptional'' case, where
$p\neq 2$ or $d\equiv 1\bmod 4$,  we have that $\nu_p(d^{ep^i}-1) =
a+i$ for all $i\geq0$,   so as a consequence of our claim, for
$s-a\leq t\leq s-1$, the order $o(p^t M_p)$ is still equal to $e$,
with $\nu_p(d^{ o(p^t M_p)}-1) = \nu_p(d^e-1)=a$, and for $t=s-a-i$
with $i\geq1$, the order $o(p^t M_p)$ is equal to $ep^i$, with
$\nu_p(d^{ o(p^t M_p)}-1) = \nu_p(d^{ep^i}-1)=a+i$. This proves the
result in the ``non-exceptional'' case.



In the ``exceptional'' case where  $p=2$ and $d\equiv 3\bmod 4$, we
have $e=1$ and $a=1$, and according to Lemma~\ref{Lpcont}, for some
integer $b\geq3$ we have that $\nu_2(d^{2+i}-1)=b+i$ for all $i\geq0$.
So as a consequence of our claim, $o(p^{s-2}M_2)=b$ and
$\nu_2(d^2-1)=b$; then for $s-b\leq t\leq s-3$, the order $o(2^t M_2)$
is still equal to $b$, with $\nu_2(d^{ o(2^t M_2)}-1) =
\nu_2(d^2-1)=b$, and for $t=s-b-i$ with $i\geq1$, the order $o(2^t
M_2)$ is equal to $e2^i$, with $\nu_2(d^{ o(2^t M_2)}-1) =
\nu_2(d^{e2^i}-1)=b+i$. This proves the result in the ``exceptional''
case.
%
%
\epf 

Now we are ready to define the $\te_v$ in the cases where
$v\in\cV^*$, that is, when $v$ is of the form $p^tM_p$ for some $p\in
\bbP$ with $0\leq t <\nu_p(m)$.  In the ``non-exceptional case'', 
i.e. for $p$ odd, or $p=2$ and $d\equiv 1\bmod 4$, and also for
$4\nmid m$, we let  
\[ \te_{M_p}=\te_{m/\pi_p(m)}=\pi_p(m)e_{M_p},\]
and  for $1\leq t<\nu_p(m)$, we let
\[ \te_{p^tM_p}=e_{p^tM_p} - \gl_{p,t} e_{M_p},\] where $\gl_{p,t}$ is
such that 
\beql{LEordcon}
\gl_{p,t}=\frac{d^{o(M_p)}-1}{d^{o(p^tM_p)}-1} \mu_{p,t} \equiv p^t
\bmod \pi_p(m)\eeql 
for some suitable integer $\mu_{p,t}$. (We will
show in a moment that this is indeed possible.) In the ``exceptional
case'' where $p=2$ and $d\equiv 3\bmod 4$ and when also $4|m$, we we
do not change the definition of~$\te_{p^tM_2}$ for $t=1, \ldots,
\pi_2(m)-2$, but we let
\beql{LEte2}\te_{M_2}=e_{2^{s-1}M_2} -2^{s-1}e_{M_2}, \qquad
\te_{2^{s-1}M_2} = \te_{m/2}= 2e_{2^{s-1}M_2},\eeql
where $s=\pi_2(m)$.  
\btm{LLte2}  
For $v\in \cV^*$,  let $\te_v $ be defined as above.
Then $\te_v\in S(m,d)$. Moreover,
$\te_v$ and $e_v$ have the same the order in~$S(m,d)$, 
except in the following two cases.

1. In the ``non-exceptional'' case where $p\in \bbP$ with $p\neq 2$ or
$d\equiv 1 \bmod 4$, and also when 
$4\nmid m$, 
the order of $\te_{M_p}=\pi_p(m)e_{M_p}$ is  
$1/\pi_p(m)$ times the order of $e_{M_p}$. 

2. In the ``exceptional''
case where $p=2$ and  $d\equiv 3 \bmod 4$, if also $4|m$ then the
order of $\te_{m/2}$ is half the order of $e_{m/2}$ and 
the order of $\te_{M_2}$ is $2/\pi_2(m)$ times
the order of $e_{M_2}$.  

Finally, the $\te_v$ for $v\in \cV$ 
are independent and generate $S(m,d)$, and Theorem~\ref{LTmain} holds.
\etm 
\bpf 
We begin by showing that the $\te_{p^tM_p}$ with $1\leq
t\leq \pi_p(v)-1$ are well-defined. To this end, first observe that by
definition of the order, we have that $o(\gl v) |o(v)$, and hence
$d^{o(\gl v)}-1$ divides $d^{o(v)}-1$. So the fraction in
(\ref{LEordcon}) is an integer, and as a consequence of
Corollary~\ref{LLordcon}, in all relevant cases the exponent of the
highest power of $p\in \bbP$ dividing this integer is at most $t$, so
that an integer $\mu_{p,t}$ for which  (\ref{LEordcon}) holds indeed
can be found. 

Next,  Theorem~\ref{Lsub} states that  $\te_{p^tM_p}-\gl_{p,t}e_{M_p}$
is in~$S(m,d)$ if and only if $p^tM_p - \gl_{p,t} M_p\equiv 0 \bmod
m$, or, equivalently, if $p^t \equiv \gl_{p,t} \bmod m$, which holds
since it is just the second requirement in (\ref{LEordcon}).  By the
same theorem, obviously $\te_{M_p} =\pi_p(m)e_{M_p}$ is also
in~$S(m,d)$, and the same holds for $\te_{M_2}$ and $\te_{m/2}$ as
defined in (\ref{LEte2}).

Then, we observe that $\te_{p^tM_p}=e_{p^tM_p}-e_{M_p}$ and
$e_{p^tM_p}$ have the same order if and only if
$(d^{o(p^tM_p)}-1)\gl_{p,t}e_{M_p}=0$, or, equivalently, if
$d^{o(M_p)}-1$ divides $(d^{o(p^tM_p)}-1)\gl_{p,t}$; this holds since it
is just the first requirement in (\ref{LEordcon}). Also, since
$\pi_p(m)$ divides the order $d^{o(M_p)}-1$ of $e_{M_p}$ for $p\in
\bbP$, we immediately have that $\pi_p(m)\te_{M_p}$ has the order as
claimed.  Now consider the exceptional case where $p=2$, $d\equiv
3\bmod 4$ and $4|m$. 
Since $e_{m/2}$ has order $d-1$, which is $2\bmod 4$,
the order of $\te_{m/2}$ is as claimed. Now let $s=\nu_2(m)$.  To
determine the order of $\te_{M_2}$, we first note that $e_{M_2}$ has
order $d^f-1$, where $f$ is the order of $d$ modulo $2^s$. Hence
$2^{s-1}e_{M_2}$ has order $(d^f-1)/2^{s-1}$. We claim that
$\te_{M_2}$ has the same order. To show this, it is sufficient to show
that the order $d-1$ of $e_{m/2}$ divides $(d^f-1)/2^{s-1}$. To this
end, note that both $(d-1)/2$ and $2^{s}$ divide $d^f-1$; since
$d\equiv 3 \bmod 4$, we know that $(d-1)/2$  is odd. Hence
$((d-1)/2, 2^s)=1$, from which the conclusion follows.

Finally, it is fairly obvious that the $\te_{p^tM_p}$ for $p\in \bbP$
and $1\leq t\leq \nu_p(m)-2$ together with $e_{M_p}$ and  $e_{m/p}$
$p\in \bbP$ are independent and generate the same subgroup as the
$e_v$ for $v\in\cV^*$. Applying Theorem~\ref{LTte1}, we
conclude that the $\te_v$ for $v\in \cV$ are independent, and thus 
form a basis of a subgroup of $S(m,d)$. From the expressions for the orders
derived above, we conclude that this subgroup has order $|\gS(m,d)|/m$;
now from Corollary~\ref{LCord} we see that in fact the $\te_v$ with
$v\in \cV$ generate $S(m,d)$. So $S(m,d)=\bigoplus_{v\in \cV} \langle
\te_v\rangle$, and now Theorem~\ref{LTmain} follows from the order
expressions.
\epf

\subsection{\label{ssect:ktzrp}Adaptations for the case of generalized Kautz
graphs}
%
Essentially, the above analysis for the case of the generalized de
Bruijn graphs only depends on the orbits of the map $x\mapsto dx$
on $\bbZ_m\setminus\{0\}$ and the order of $d$ modulo powers of $p$.
Thus, similar results hold when $d<0$; the only difference is that
relevant group elements $g$ now have orders of the form
$|d^e-1|/c$, so that the group $\langle g \rangle$ is 
isomorphic to $\bbZ_{|d^e-1|/c}$. 
We leave the details for the interested reader.


\section{The group of invertible circulant matrices}

In \cite{DuPa-neck} it was shown that a family of bijections between 
the set of 
aperiodic necklaces of length~$n$ over the finite field $\bbF_q$, 
with $q$ prime, and the set of degree $n$ normal polynomials over $\bbF_q$,
gives rise to a permutation group on any of these sets.
This group turns out
to be isomorphic to $C(n,q)$; as well, \cite{DuPa-neck} conjectured
a relation of $C(n,q)$ and $S(n,q)$. 
The present work confirms this relation\footnote{In view of this
relation, it would be desirable to have an explicit embedding of the
group $S(n,q)$ as a subgroup of~$C(n,q)$.} 
with the sandpile group $S(n,q)$ of $\mathrm{DB}(n,q)$.
We show that indeed $C'(n,q)/\langle x\rangle$ is isomorphic
to~$S(n,q)$ for all $n$ and all primes $q$, and also for all prime
powers $q$ provided that $(n,q)=1$, by explicitly computing a
decomposition into cyclic subgroups. 
While the case $(n,q)=1$ appears to be well-understood, we did not
find an explict reference in the literature.
The general case is harder, and the only relevant reference we found was
the case $n=2^k$, $p=2$ dealt with in \cite[Prop. XI.5.7]{MR0249491}.

In what follows, we use the description for $C'(n,q)$ as given in
(\ref{LECp}). 
\btm{LTiso-reg-circ} Let $q$ be a prime power, and let $m$ be a
positive integer with $(m,q)=1$. 
Then \[C'(m,q)\cong \gS(m,q)\] and \beql{LE-reg-circ} C'(m,q) /
\langle x \rangle \cong S(m,q).  \eeql
\etm 

\bex{LEnotgen}
We remark that the condition
$(m,q)=1$ (respectively $(n,d)=1$) is necessary in 
Theorem~\ref{LTiso-reg-circ} 
(respectvely in Theorem~\ref{LTmain-circ}).
Indeed, one may check directly that
$C'(9,9)/\langle x\rangle\cong \mathbb{Z}_9^4\oplus
\mathbb{Z}_{27}\oplus\mathbb{Z}_3^3$, 
while $S(9,9)\cong \mathbb{Z}_9^7$.
\eex

\bpf Let $\cV$ be defined as in the beginning of
Section~\ref{sect:rp}, so that $\cV^+=\cV\cup\{0\}$  is
a complete set of orbit representatives for the map $x\mapsto qx$
on $\bZ_m$ containing the set $\cV^*$ of all numbers of the form
$p^iM_p$ for a prime $p\in \bbP$, with $M_p=m/\pi_p(m)$. (Refer to the
beginning of Section~\ref{sect:rp} for definitions of
$\cV^*$, $M_p$, and $\pi_p(m)$.)
We write $o(v)$ to denote the size of the orbit of $v\in\cV^+$. 
Next, let $q$ have order $k$ modulo $q$, so that $n|q^k-1$, and let
$\ga$ be a primitive element of~$\GF(q^k)$; put $\xi=\ga^{(q^k-1)/n}$.
For $v\in \cV^+$, let $f_v(x)$ denote the minimal polynomial of
$\xi^{v}$ over $\bbF_q$; note that $f_v(x)=\prod_{i=0}^{o(v)-1}
(x-\xi^{q^iv})$ is $\bbF_q$-irreducible of degree~$o(v)$. 
Then we can write $x^m-1=\prod_{v\in \cV^+}
f_v(x)$. Thus
$C(m,q)\cong \displaystyle\prod_{v\in \cV^+} 
\left(\bbF_q[x]/( f_v(x))\right)^*$ and 
\[C'(m,q)\cong \prod_{v\in \cV}\left(\bbF_q[x] /(f_v(x))\right)^* 
    \cong \oplus_{v\in \cV} \bbF_{q^{o(v)}}^* 
   \cong \oplus_{v\in \cV} \bbZ_{q^{o(v)}-1}.\] 
Hence $C'(m,q)\cong \gS(m,q)$ by Theorem~\ref{TSd}.

Next, we consider the quotient $C'(m,q)/\langle x \rangle$. The image
of $x\bmod f_v(x)$ in~$\bbF_{q^{o(v)}}$ is $\xi^{v}$; since $\xi$ has
order~$m$, we see that $x$ has order~$m/(m,v)$ in~$\bbF_{q^{o(v)}}$.
Hence for each $v\in \cV$,
we can choose a primitive element $\gb_v$
in~$\bbF_{q^{o(v)}}\cong\bbF_q[x]\bmod f_v(x)$, so with $\gb_v$ of
order $q^{o(v)}-1$,  such that the image of $x$  in~$\bbF_{q^{o(v)}}$
is $\gb^{r_v}$, where \beql{LEfi}r_v=\frac{q^{o(v)}-1}{m/(m,v)};\eeql
as a consequence, we have that
\beql{LE-sp-cong} G:=C'(m,q) / \langle x
\rangle \cong G= \langle \gb_v \mid v\in \cV, 
\mbox{$\ord(\gb_v)=q^{o(v)}-1$ ($v\in \cV$) and $\prod_{v\in \cV}
\gb_v^{r_v}=1$} \rangle.  
\eeql

To obtain the group structure  of $G$ in (\ref{LE-sp-cong}), we 
investigate the Sylow-$p$ subgroup $G_p\leq G$ for each prime $p$.
After some standard manipulations, i.e. 
writing $\gb_v=\prod_p  \gb_{v,p}$ with
$\ord(\gb_{v,p})=\pi_p( q^{o(v)}-1)$, then fixing a prime $p$ and
letting $\gc_v=\gb_{v,p}$, 
we obtain that \[G_p = \langle \gc_v\mid v\in \cV, 
\mbox{$\ord(\gc_v)=\pi_p(q^{o(v)}-1)$ ($v\in\cV$) and $\prod_{v\in
\cV} \gc_v^{s_v}=1$} \rangle,\] where 
$s_v=\pi_p(r_v)=\pi_p((q^{o(v)}-1)(v,m)/m)$ for all $v\in \cV$. First,
observe that if \mbox{$p\!\!\not | m$}, then
$s_v=\pi_p(q^{o(v)}-1)=\ord(\gc_v)$, and hence \[G_p= \bigoplus_{v\in
\cV} \bbZ_{\pi_q(p^{o(v)}-1)}.\]
Now, let $p|m$. We need to determine for which $v\in
\cV$ the number $s_v$ is minimal. 



\mymedskip

\noindent 
{\bf Claim 1}: If $p|m$ and $\pi_p(v)=p^k$, then $s_v\geq s_{p^kM_p}$.

\mymedskip


\noindent  
Indeed, first note that $p^kM_p(p^e-1)\equiv 0\bmod m$ precisely when
$p^e-1\equiv 0 \bmod \pi_p(m)/p^k$.
Then, writing $v=p^kw$ with $(p,w)=1$, we see that $p^kw(p^e-1)\equiv
0 \bmod m$ implies that 
$(p^e-1) \equiv 0 \bmod \pi_p(m)/p^k$, and the claim is now obvious.

\mymedskip

\noindent {\bf Claim 2}: If $p\neq 2$ or $q\not \equiv 3 \bmod 4$ or
$4\!\!\not | m$, then $s_v$ is minimal for $v=M_p=m/\pi_p(m)$. 

\mymedskip

\noindent  Indeed, if $p=2$ and \mbox{$4\!\!\not | m$}, the claim
follows immediately from claim 1.  Otherwise, let $v_{j}=m/p^j$, and
let  $\gl_j=\pi_p(q^{o(v_j)}-1)$. 
As a consequence of  Lemma~\ref{Lpcont}, we have that \[(\gl_1, \gl_2,
\ldots ) = (q^a, \ldots, q^a, q^{a+1}, q^{a+2}, \ldots),\] where $q^a$
occurs $a$ times, for some $a\geq 1$. So, since $s_{v_j}=\gl_j
p^{-j}$, we see that $s_{v_j}$ is minimal when $j$ is as large as
possible. Now the claim follows from claim 1.

To complete the determination of~$G_p$ in the non-exceptional case where
$p\neq 2$ or $q\not \equiv 3 \bmod 4$, or when $4\nmid m$, we
proceed as follows. 
Define $\gd = \gc_{M_p} \prod_{v\neq M_p} \gc_v^{s_v/s_{M_p}}$.  Note
that $\gd$ has order $s_{M_p}=\pi_p(q^{M_p}-1)/\pi_p(m)$ in $G_p$.
Obviously, $\gd$ and the $\gc_v$ with $v\neq M_p$ generate $G_p$, so
we conclude that 
\[G_p\leq
\bbZ_{\pi_p(q^{o(M_p)}-1)/\pi_p(m)} \oplus \bigoplus_{v\in
\cV\setminus\{M_p\}} \bbZ_{\pi_p(q^{o(v)}-1)}.\]



Now consider the exceptional case where $p= 2$ and $q \equiv 3 \bmod
4$ and $4|m$.  With the same notation as after claim 2,
Lemma~\ref{Lpcont} now implies that \[(\gl_1, \gl_2,\ldots) = (2, 2^b,
\ldots, 2^b, 2^{b+1}, 2^{b+2}, \ldots),\] where $2^b$ occurs $b-1$
times, for some $b\geq 3$. So \[(s_{v_1}, s_{v_2}, ...) = (1, 2^{b-2},
\ldots, 2, 1, 1, \ldots).\] Since $v_1=m/2$, we see that $s_{m/2}=1$
and $o(m/2)=1$, so $\ord(m/2)=\pi_2(q^{o(m/2)}-1)=\pi_2(q-1)=2$. Also, 
note that $s_{M_2}= \pi_2(q^{o(M_2)}-1)/\pi_2(m)$. Now define
$\gd=\gc_{M_2}\prod_{v\neq M_2, m/2} \gc_v^{s_v/s_{M_2}}$. Then
$\gd^{s_{M_2}}=\gc_{m/2}^{-s_{m/2}}=\gc_{m/2}$, and since
$\ord(\gc_{m/2})=2$, we conclude that $\gd$ has order $2s_{M_2}$.
Also, since $s_{m/2}=1$, $\gc_{m/2}$ can be expressed in terms of
elements $\gc_v$ with $v\neq m/2$. It is now easy to see that $G_p$ is
generated by $\gd$ and $\gc_v$ for $v\in\cV \setminus\{M_2,  m/2\}$,
and since $\pi_2(q^{o(m/2)}-1)=2$, we conclude that in this case
\[G_2\leq\bbZ_{\pi_p(q^{o(M_2)}-1)/\pi_2(m)}
\oplus \bbZ_{(q^{o(m/2)}-1)/2} \oplus \bigoplus_{v\in
\cV\setminus\{M_2,m/2\}} \bbZ_{\pi_2(q^{o(v)}-1)}.\]

If we combine the information about the various Sylow $p$-subgroups,
we may conclude from Theorem~\ref{LTmain} that $G$ is a subgroup of
$S(m,q)$. Since $|G|=|C'(m,q)/\langle x \rangle|=
|C'(m.q)|/m=|S(m,q)|$, we conclude that $G=S(m,q)$ as desired.  \epf

\section{The group $C(n,p)$ for general $n$ and prime $p$}

\newcommand{\ZZ}{\bbZ}
We work with $C(n,p)$ as the group of invertible circulant matrices
over $\FF_p$ and with $C^*(n,p)\leq C(n,p)$ as the subgroup of marices
with eigenvalue $1$ on the eigenvector $(1, \ldots, 1)$.

We first compute the group decomposition of $C^*(n,p)$ and then
proceed to compute the group decomposition of $C^*(n,p) / \langle Q_n
\rangle $.  

\ble{LLsh1} Let $p$ be prime and $(m,p)=1$. Then 
\[C^*(p^km,p)= \left[
\bigoplus _{i=0}^{k-2} \mathbb{Z}_{p^{k-1-i}}^{p^{i} (p-1)^2 m}
\right]  \oplus \ZZ_{p^k}^{(p-1)m} \oplus C^*(p,m).\]
\ele

\bpf Note that  $\phi : C^*(pn,p) \to C^*(pn,p)$ given by $x \mapsto
x^p$ is a well-defined homomorphism, as $\phi$ preserves both the
determinant and the eigenvalue of the common eigenvector
$(1, \ldots 1)$.  It can be easily checked that
$|\ker(\phi)|=p^{(p-1)n}$. 

It can be also easily checked that $(C^*(pn,p))^p$ is isomorphic to 
$C^*(n,p)$, via $\sum_{i=0}^{n-1} a_i Q_{pn}^{p i}  \mapsto
\sum_{i=0}^{n-1} a_i Q_n^{i} .$ Hence $\phi$ can be viewed as
surjective homomorphism from $C^*(pn,p)$ to $C^*(n,p)$. 

We prove the lemma by using induction on $k$. It is easy to see
that $C^*(pm,p)=\ZZ_p^{a_1} \oplus C^*(m,p)$ since $(C^*(pm,p))^p
\simeq C^*(m,p)$.  Since in this case $|\ker(\phi)|=p^{(p-1)m}$, we
can conclude that $a_i=(p-1)m$ and the statement is true for
$k=1$. 

By using the fact that $(C^*(p^{k+1}m,p))^p \simeq C^*(p^km,p)$, by
induction we conclude that 
\[
C^*(p^{k+1}m,p) = Z_p^{a_1} \oplus
\left[ \bigoplus _{i=0}^{k-2} \mathbb{Z}_{p^{k-i}}^{p^{i} (p-1)^2 m}
\right]  \oplus \ZZ_{p^{k+1}}^{(p-1)m} \oplus C^*(p,m).
\]

It remains to find $a_1$.  Note that when
$n=p^{k+1}m$, we have 
\[
\ker (\phi)= Z_p^{a_1} \oplus \left[ \bigoplus
_{i=0}^{k-2} \mathbb{Z}_{p}^{p^{i} (p-1)^2 m} \right]  \oplus
\ZZ_{p}^{(p-1)m}.
\]
Recall that 
$|\ker(\phi)|=p^{(p-1)p^km}$. Thus we have two expressions for 
$|\ker(\phi)|$.  Equating them, we have: 
\begin{align*}
(p-1)p^km&= a_1 + \left[ \sum_{i=0}^{k-2}p^i (p-1)^2m \right] +(p-1)m,
\\ (p-1)p^km &= a_1 + (p^{k-1}-1)(p-1)m +(p-1)m, \\ 
 a_1&= p^{k-1}(p-1)^2m, 
\end{align*}
which completes the induction.
\epf

Finally, we complete the proof of Theorem~\ref{LTmain-circ}.
Let $d=p$ be prime, $n=p^km$ and $(p,m)=1$.
Since $\langle Q_n \rangle$ is cyclic of order $p^km$, 
its Sylow $p$-subgroup is $Z_{p^k}$. 
Combined with Lemma~\ref{LLsh1}, this implies that 
the Sylow $p$-subgroup of the quotient group $C^*(n,p)/\langle
Q_n \rangle$ equals
\[ \mathrm{Syl}_p( C^*(p^{k}m,p)/ \langle Q_n
\rangle)  = \left[ \bigoplus _{i=0}^{k-2}
\mathbb{Z}_{p^{k-1-i}}^{p^{i} (p-1)^2 m} \right]  \oplus
\ZZ_{p^{k}}^{(p-1)m-1}.\]

By Lemma~\ref{LTcrt} and Theorem~\ref{LTmain1} with $d=p$, 
the Sylow $p$-subgroup of $S(n,p)$ satisfies
\[\mathrm{Syl}_p(S(n,p))\cong S_0(n,p)
 \cong \bigoplus_{i=0}^{k-1}\biggl(\bbZ_{p^{i+1}/d_i}\oplus
\bbZ_{p^{i+1}}^{n_i-2n_{i+1}+n_{i+2}-1}\biggr),\] 
with $n=n_0,n_1,\dots,n_k=n_{k+1}$ the $p$-sequence of $n$, 
and $d_i=n_i/n_{i+1}=(n_i,p)$ for $0\leq i\leq k$.
As $n=p^km$, we have $n_i=p^{k-i}m$ and $d_i=p$ for $0\leq i\leq k$, 
and $n_i=m$ otherwise.
Plugging these values into above equation, we have
\[
\mathrm{Syl}_p(S(p^km,p))=  \ZZ_{p^{k-1}} \oplus \ZZ_{p^k }^{(p-1)m-1}
\oplus \bigoplus_{i=0}^{k-2}   \left[ \ZZ_{p^{i}} \oplus
\ZZ_{p^{i+1}}^{p^{k-2-i}(p-1)^2m-1} \right].
\]
Changing $i$ to $k-2-i$ in the above summation, we have
\begin{align*}
\mathrm{Syl}_p(S(p^km,p))&=  \ZZ_{p^{k-1}} \oplus \ZZ_{p^k }^{(p-1)m-1}
\oplus \bigoplus_{i=0}^{k-2}   \left[ \ZZ_{p^{k-2-i}} \oplus
\ZZ_{p^{k-1-i}}^{p^{i}(p-1)^2m-1} \right]\\ &= \left[ \bigoplus
_{i=0}^{k-2} \mathbb{Z}_{p^{k-1-i}}^{p^{i} (p-1)^2 m} \right]  \oplus
\ZZ_{p^{k}}^{(p-1)m-1}\\ 
&= \mathrm{Syl}_p(C^*(p^{k}m,p)/ \langle Q_n\rangle)
\end{align*} 
and the proof of Theorem~\ref{LTmain-circ} is complete.

\bibliographystyle{abbrv} 
\bibliography{sandpile-euro}

\end{document}